\documentclass[reqno]{amsart}
\usepackage{graphicx,amssymb}

\newcommand{\FF} {{\mathbb F}}

\newcommand{\QQ} {{\mathbb Q}}
\newcommand{\NN} {{\mathbb N}}
\newcommand{\RR} {{\mathbb R}}
\newcommand{\PP} {{\mathbb P}}
\newcommand{\ZZ} {{\mathbb Z}}

\newcommand{\m}{{\mathfrak m}}

\newcommand{\xx} {{\mathbf x}}
\newcommand{\yy} {{\mathbf y}}

\newcommand{\isthmus} {{\mathrm{isthmus}}}
\newcommand{\looop} {{\mathrm{loop}}}
\newcommand{\rank}{{\mathrm{rank}}}

\newcommand{\std}{{\mathrm{std}}}
\newcommand{\Des}{{\mathrm{Des}}}

\newcommand{\LL}{{\mathcal L}}
\newcommand{\AAA}{{\mathcal A}}
\newcommand{\Mat}{{\mathcal Mat }}
\newcommand{\Qsym}{{\mathcal{QS}ym }}
\newcommand{\bases}{{\mathcal B}}
\newcommand{\flag}{{\mathcal F}}

\newcommand{\Fdual}{{{F}^{\ast}}}
\newcommand{\fdual}{{{\phi}^{\ast}}}

\newtheorem{theorem}{Theorem}
\newtheorem{corollary}[theorem]{Corollary}
\newtheorem{proposition}[theorem]{Proposition}

\newtheorem{lemma}[theorem]{Lemma}
\newtheorem{defn}[theorem]{Definition}
\newtheorem{example}[theorem]{Example}
\newtheorem{remark}[theorem]{Remark}
\newtheorem{question}[theorem]{Question}
\newtheorem{problem}[theorem]{Problem}

\numberwithin{theorem}{section}
\numberwithin{equation}{section}

\begin{document}
\title{A quasisymmetric function for matroids}

\author{Louis J. Billera}
\address{Department of Mathematics\\
Cornell University\\
Ithaca NY 14853-4201}
\email{billera@math.cornell.edu}

\author{Ning Jia}
\address{Department of Mathematics\\
Virginia Tech\\
460 MCB\\
Blacksburg, VA 24060}
\email{njia@vt.edu}

\author{Victor Reiner}
\address{School of Mathematics\\
University of Minnesota\\
Minneapolis, MN 55455}
\email{reiner@math.umn.edu}


\thanks{First author supported by NSF grant DMS-0100323; third author supported
by NSF grant DMS-0245379.
This work was begun while the first author was an Ordway Visitor
at the University of Minnesota School of Mathematics.
The authors thank Marcelo Aguiar, Saul Blanco, Harm Derksen, Sam Hsiao, Sean Keel, Kurt Luoto, Frank Sottile, David Speyer, and
Jenia Tevelev for helpful conversations, and are grateful to Alexander Barvinok, Jim Lawrence 
and Peter McMullen for helpful suggestions on Proposition \ref{polar_relation}.}

\begin{abstract}
A new isomorphism invariant of matroids is introduced, in the form of a quasisymmetric function.
This invariant 
\begin{enumerate}
\item[$\bullet$]
defines a Hopf morphism from the Hopf algebra of matroids to the quasisymmetric functions,
which is surjective if one uses rational coefficients,
\item[$\bullet$]
is a multivariate generating function for integer weight vectors
that give minimum total weight to a unique base of the matroid,
\item[$\bullet$]
is equivalent, via the Hopf antipode, to a generating function for integer weight vectors
which keeps track of how many bases minimize the total weight,
\item[$\bullet$]
behaves simply under matroid duality,
\item[$\bullet$] 
has a simple expansion in terms of $P$-partition enumerators, and
\item[$\bullet$]
is a valuation on decompositions of matroid base polytopes.
\end{enumerate}
This last property leads to an interesting application: 
it can sometimes be used to prove that a matroid
base polytope has no decompositions into smaller matroid base polytopes.   Existence of
such decompositions is a subtle issue arising in work of Lafforgue,
where lack of such a decomposition implies the matroid has only a finite number of realizations
up to scalings of vectors and overall change-of-basis.
\end{abstract}

\maketitle

\tableofcontents

\section{Definition as generating function}
\label{definition-section}
We begin by defining the new matroid invariant.  For matroid terminology 
undefined here, we refer the reader to some of the standard references, such
as \cite{CrapoRota, Oxley, Welsh, White1, White2, White3}.

Let $M=(E,\bases)$ be a matroid on ground set $E$, with bases $\bases=\bases(M)$.
Let $\PP:=\{1,2,3,\ldots\}$ be the positive integers. 
We will say that a weighting function $f: E \rightarrow \PP$ is {\it $M$-generic} if
the minimum $f$-weight $f(B):=\sum_{e \in B} f(e)$ among all bases $B$ of $M$
is achieved by a {\it unique} base $B \in \bases(M)$.  For example, it is a standard
exercise in matroid theory (see, e.g. \cite[Exer. 1.8.4]{Oxley}) to show that 
$f$ is $M$-generic if $f$ is {\it injective}, that is, if $f$ assigns all distinct weights.

\begin{defn} \rm \ 
\label{definition-of-F}
Given a matroid $M$ as above,
define a power series $F(M,\xx)$ in countably many variables
$x_1,x_2,\ldots$ as the generating function for $M$-generic weighting functions $f$
according to number of times $f$ takes on each value in $\PP$.  That is,
\begin{equation}
\label{definition-equation}
F(M,\xx):=\sum_{\substack{ M{\textrm{-generic }} \\ f: E \rightarrow \PP}} \xx_{f}
\end{equation}
where $\xx_{f}:=\prod_{e \in E} x_{f(e)}$.
\end{defn}

One of the defining properties of a matroid \cite[Theorem 1.8.5]{Oxley} is that
an $f$-minimizing base may be found by {\it (Kruskal's) greedy algorithm}:
\begin{quote}
Construct a sequence of independent sets 
$$
\emptyset = :I_0, I_1, \ldots, I_{\rank(M)}
$$
by defining $I_j:=I_{j-1} \cup \{e\}$ where $e$ is any
element in $E$ having minimum weight $f(e)$ among those for which $I_{j-1} \cup \{e\}$
is independent.  Then $I_{\rank(M)}$ is an $f$-minimizing base of $M$.
\end{quote}

\section{Quasisymmetry}
\label{quasisymmetry-section}

We recall \cite{Gessel}, \cite[\S 7.19]{Stanley-EC2} 
what it means for a power series $f(\xx)$ in a linearly ordered variable set
$x_1,x_2,\ldots$ to be {\it quasisymmetric}:  $f$ must have bounded degree, and for any fixed composition 
$(\alpha_1,\ldots,\alpha_k)$ in $\PP^k$,
the coefficient of the monomials 
$x_{i_1}^{\alpha_1} x_{i_2}^{\alpha_2} \cdots x_{i_k}^{\alpha_k}$ with $i_1 < i_2 < \cdots < i_k$ 
in $f$ are all the same.  
Put differently, $f$ is quasisymmetric if and only if it is a (finite) linear combination of the
{\it monomial} quasisymmetric functions\footnote{While there is a danger of confusion between
matroids $M$ and monomial quasisymmetric functions $M_{\alpha}$, the difference will
always be clear by the context.} indexed by compositions
$\alpha=(\alpha_1,\ldots,\alpha_k)$:
$$
M_\alpha:=\sum_{ 1 \leq i_1 < i_2 < \cdots < i_k} 
 x_{i_1}^{\alpha_1} x_{i_2}^{\alpha_2} \cdots x_{i_k}^{\alpha_k}.
$$

\begin{proposition}
\label{quasisymmetry}
For any matroid $M$, the power series $F(M,\xx)$ is {\it quasisymmetric}.
\end{proposition}

\begin{proof}
This follows from the fact that the $f$-minimum bases can all be found by the
greedy algorithm, and this algorithm makes all of its decisions based only on
the {\it relative ordering} and equality of various weights $f(e)$, not on their actual values.
\end{proof}

\begin{example} \rm \ 
\label{tiny-examples}
When $|E|=0$, there is only one matroid $M_\emptyset$, having rank $0$ and exactly one base,
the empty base $\emptyset$.  As there is only one function $f$ from the empty set $E$ into
$\PP$, and this $f$ has no coordinates (!), we should decree
$\xx_f=1$ (as the empty product is $1$). Hence $F(M_\emptyset,\xx)=1$.

There are two matroids with $|E|=1$, namely $M_\isthmus$ of rank $1$ having a single base
$\{e\}$, and $M_\looop$ of rank $0$ having a single base $\emptyset$.
Every $f:E \rightarrow \PP$ is generic for either of these, so that
$$
F(M_\isthmus,\xx) = F(M_\looop,\xx) = x_1 + x_2 + x_3 + \cdots = M_1.
$$
\end{example}

The enumerative information recorded in $F(M,\xx)$ is data about optimizing
weight functions on the bases of $M$.  An obvious specialization counts $M$-generic
weight functions that take on only a limited number of distinct weight values.

\begin{defn} \rm \
For a positive integer $m$, let $[m]:=\{1,2,\ldots,m\}$, and define
$$
\begin{aligned}
\phi(M,m)
      &: =  F(M, \xx)_{\substack{x_1 = x_2 =\cdots = x_m =1 \\ x_{m+1}=x_{m+2}= \cdots =0}} \\
      & = |\{ M\textrm{-generic }f:E \rightarrow [m] \}|.
\end{aligned}
$$
\end{defn}
Since $F(M,\xx)$ is a power series of bounded degree, $\phi(M,m)$ is a polynomial function of $m$.
When $m$ is large, almost all weight functions $f:E \rightarrow [m]$ are injective and hence 
$M$-generic, so the polynomial expansion of $\phi(M,m)$ begins
$$
\phi(M,m) = m^n + O(m^{n-1})
$$
where $n:=|E|$.

  In Section~\ref{antipodes}, our analysis of the behavior of $F(M,\xx)$ under
the Hopf algebra antipode on quasisymmetric functions will
imply an interesting reciprocity result for the polynomial $\phi(M,m)$.

  In the remainder of this paper, we will suppress the $\xx$ in $F(M,\xx)$ and write $F(M)$
unless there is a need to consider the variables.

\section{Hopf algebra morphism}
\label{Hopf-section}

  There is a known Hopf algebra structure built from matroids 
\cite{CrapoSchmitt1, CrapoSchmitt2, CrapoSchmitt3, Schmitt}  
and a perhaps better-known Hopf algebra of quasisymmetric functions \cite[\S 4]{Gessel}.   
The goal of this section is to show that the invariant $F(M)$ defines a Hopf morphism
between them.

  Let $\Mat$ be the free $\ZZ$-module  
consisting of formal $\ZZ$-linear combinations
of basis elements $[M]$ indexed by isomorphism classes of matroids $M$.  Endow $\Mat$ with
a product and coproduct extended $\ZZ$-linearly from the following definitions on
basis elements:

\begin{equation*}
\begin{aligned}
\, [ M_1 ] \cdot [ M_2 ] & := [M_1 \oplus M_2] \\
\Delta [M] :& = \sum_{A \subseteq E} [M|_A] \otimes [M/A]
\end{aligned}
\end{equation*}

\noindent
where $M_1 \oplus M_2$ is the {\it direct sum} of the matroids $M_1, M_2$, and
$M|_A, M/A$ denote the {\it restriction} of 
$M$ to $A$ and the {\it contraction} (or {\it quotient}) of $M$ by $A$,
respectively.  One has a $\ZZ$-module direct sum decomposition
$\Mat = \bigoplus_{n \geq 0} \Mat_n$, where $\Mat_n$ denotes the submodule spanned
by the basis elements $[M]$ for which the ground set $E$ of $M$ 
has cardinality $|E|=n$.  One can then 
easily check that this product and coproduct make 
$\Mat$ into a graded, connected Hopf algebra over $\ZZ$ which is commutative, but non-cocommutative.
Here the unit is $[M_\emptyset]$.

  Let $\Qsym$ (or $\Qsym({\bf x})$)
denote the Hopf algebra of quasisymmetric functions in the linearly ordered
variable set $x_1,x_2,\ldots$ and having coefficients in $\ZZ$.  The product in $\Qsym$ is 
inherited from the formal power series ring $\ZZ[[x_1,x_2,\ldots]]$.  The coproduct may be
described as follows.  A quasisymmetric function $f(\xx)$ defines a unique quasisymmetric
function $f(\xx,\yy)$ in the linearly ordered variable set 
$$
x_1 < x_2 < \cdots < y_1 < y_2 < \cdots
$$
by insisting that $f(\xx, {\mathbf 0})=f(\xx)$.  In other words, for any $i_1 < \cdots < i_k$ and
$j_1 < \cdots < j_\ell$, the coefficient of 
$x_{i_1}^{\alpha_1} \cdots x_{i_k}^{\alpha_k} y_{j_1}^{\beta_1} \cdots y_{j_\ell}^{\beta_\ell}$
in $f(\xx,\yy)$ is defined to be the coefficient of 
$x_1^{\alpha_1} \cdots x_k^{\alpha_k} x_{k+1}^{\beta_1} \cdots x_{k+\ell}^{\beta_\ell}$
in $f(\xx)$.  Consider the injective map
$$
i: \ZZ[[x_1,x_2,\ldots]] \otimes \ZZ[[y_1,y_2,\ldots]] \rightarrow \ZZ[[x_1,x_2,\ldots,y_1,y_2,\ldots]] 
$$
which sends $f(\xx) \otimes g(\yy)$ to $f(\xx)g(\yy)$.  The image $i(\Qsym({\bf x}) \otimes \Qsym({\bf y}) )$
contains the quasisymmetric functions $\Qsym({\bf x, y})$, that is, there
is a unique expansion $f(\xx,\yy) = \sum_i f_i(\xx) g_i(\yy)$ for any quasisymmetric function $f(\xx,\yy)$.
This defines the coproduct $\Delta: \Qsym \rightarrow \Qsym \otimes \Qsym$.
Grading $\Qsym$ by the usual notion of degree, one can check that
$\Qsym$ becomes a graded, connected Hopf algebra over $\ZZ$ which is commutative, but non-cocommutative.

\begin{theorem}
\label{morphism-theorem}
The map $F$
$$ \begin{matrix}
\Mat &\rightarrow &\Qsym \\
[M] &\mapsto &F(M)
\end{matrix}
$$
is a morphism of Hopf algebras.
\end{theorem}
\begin{proof}
Example~\ref{tiny-examples} shows that $F$ sends the unit $M_\emptyset$ of $\Mat$
to the unit $1$ of $\Qsym$. The fact that $F$ preserves degree
shows that it preserves the counit.

The fact that $F$ preserves the product structures follows because the bases of
$M_1 \oplus M_2$ are the disjoint unions $B_1 \sqcup B_2$ of a base $B_1, B_2$
from each.  This implies that $f:E_1 \sqcup E_2 \rightarrow \PP$ is 
$(M_1 \oplus M_2)$-generic if and only if $f|_{E_i}$ is $M_i$-generic for $i=1,2$.

The fact that $F$ preserves the coalgebra structure is somewhat more interesting.
Unravelling the definitions, this amounts to checking the following identity:
\begin{equation}
\label{coalgebra-map}
F(M, \xx, \yy) = \sum_{ A \subseteq E} F(M|_A, \xx) \, F(M/A, \yy).
\end{equation}
The left side of \eqref{coalgebra-map} has the following interpretation.  Linearly
order the disjoint union $\PP \sqcup \PP$ as follows:
$$
1 < 2 < 3 < \cdots < 1' < 2' < 3' <\cdots 
$$
Given a weight function $f: E \rightarrow \PP \sqcup \PP$, define
$(\xx\yy)_f:=\prod_{e \in E} z_e$ where
$$
z_e := \begin{cases}
x_i & \textrm{ if }f(e) = i \textrm{ ( with no prime ) } \\
y_i & \textrm{ if }f(e) = i'
\end{cases}. 
$$
Then 
$$
F(M, \xx, \yy) = \sum_{\substack {M\textrm{-generic}\\f: E \rightarrow \PP \sqcup \PP}} (\xx\yy)_f.
$$
On the other hand, the right side of \eqref{coalgebra-map} expands to  
$\sum_{(A,f_1,f_2)} \xx_{f_1} \yy_{f_2}$, where the sum ranges over all triples 
$(A,f_1,f_2)$ in which 
\begin{enumerate}
\item[$\bullet$] $A$ is a subset of $E$,
\item[$\bullet$] $f_1: A \rightarrow \PP$ is $M|_A$-generic, and
\item[$\bullet$] $f_2: E \backslash A \rightarrow \PP$ is $M/A$-generic.
\end{enumerate}
There is an obvious association $f \mapsto (A,f_1,f_2)$ defined by
$$
\begin{aligned}
A &:= \{e \in E: f(e) \textrm{ has no prime }\} \\
f_1 &:= f|_A \\
f_2 &:= f|_{E \backslash A}.
\end{aligned}
$$
It only remains to check that $f$ is $M$-generic if and only if 
$f|_A$ and $f|_{E \backslash A}$ are $M|_A$ and $M/A$-generic, respectively.  This follows from the
sequential nature of the greedy algorithm:  because the primed values $i'$ are bigger than all
the unprimed values $i$, when the greedy algorithm finds $f$-minimizing bases for $M$,
it must first find $f|_A$-minimizing bases for $M|_A$ by trying to use only $e$'s
with unprimed values for as long as it can, and then proceed to 
find $f|_{E\backslash A}$-minimizing
bases for $M/A$ using primed values.   Lack of uniqueness in the $f$-minimizing bases of $M$
can only occur if it occurs in one of these two steps, leading either to lack of
uniqueness in the $f|_A$-minimizing bases of $M|_A$ or in the $f|_{E \backslash A}$-minimizing
bases of $M/A$.  Conversely, lack of uniqueness in either step will lead to
lack of uniqueness for the whole computation.
\end{proof}

  It turns out that the Hopf morphism $\Mat \rightarrow \Qsym$ is
{\it not} surjective if one works over $\ZZ$, but becomes
surjective after tensoring with the rationals.   The somewhat technical proof
of this surjectivity\footnote{A shortening of parts of this proof has been found
recently by Luoto, as an application of his ``matroid-friendly'' basis of quasisymmetric
functions; see \cite[\S 7.4]{Luoto}.}
is given in the Appendix (Section~\ref{appendix}).  The proof involves the
construction of two new $\ZZ$-bases for $\Qsym$, which may be of independent
interest.

\begin{remark}   ({\it on combinatorial Hopf algebras}) \rm \ \\ 
  Definition~\ref{definition-of-F} for $F(M)$ immediately implies that for any
composition $\alpha = (\alpha_1,\ldots,\alpha_k)$ of $n:=|E|$, the 
coefficient $c_\alpha$ in the unique expansion 
\begin{equation}
\label{M-expansion}
F(M)=\sum_{\alpha} c_\alpha M_\alpha
\end{equation}
has the following interpretation:  $c_\alpha$ is the number of $M$-generic
$f:E \rightarrow \PP$ in which $|f^{-1}(i)|=\alpha_i$ for $i=1,2,\ldots,k$.

The work of Aguiar, Bergeron and Sottile \cite{AguiarBergeronSottile}
on combinatorial Hopf algebras also offers an interpretation for $c_\alpha$, using the fact that
$F$ is a Hopf morphism, as we explain here.  In their theory, the {\it character} (= multiplicative linear
functional) $\zeta_{\mathcal Q}:\Qsym \rightarrow \ZZ$ 
defined by
$$
\zeta_{\mathcal Q}(M_\alpha) = 
\begin{cases} 
  1 & \text{ if }\alpha\text{ has at most one part, and} \\
  0 & \text{ otherwise}
\end{cases}
$$
plays a crucial role, making $\Qsym$ into what they call a {\it combinatorial Hopf algebra}.
The Hopf morphism $F:\Mat \rightarrow \Qsym$ then allows one to uniquely
define a character  $\zeta_{\mathcal M}:\Mat \rightarrow \ZZ$, via
$\zeta_{\mathcal M} := \zeta_{\mathcal Q} \circ F$, so that $F$ becomes a {\it morphism of
combinatorial Hopf algebras}.  

It is not hard to see directly 
(or one can appeal to Corollary~\ref{zeta-computation} below) the following more explicit
description of the character $\zeta_{\mathcal M}$.  Say that a matroid $M$ {\it splits completely} if it
is a direct sum of matroids on $1$ element, that is, a direct sum of loops and isthmuses,
or equivalently, if it has only one base $B$.
Then for any matroid $M$
$$
\zeta_{\mathcal Q}([M]) = 
\begin{cases} 
  1 & \text{ if }$M$\text{ splits completely, and} \\
  0 & \text{ otherwise.}
\end{cases}
$$
Using this, \cite[Theorem 4.1]{AguiarBergeronSottile} immediately implies another
interpretation for the coefficient $c_\alpha$ in \eqref{M-expansion}.  Given
a flag $\flag$ of subsets 
\begin{equation}
\label{typical-flag}
\flag: \emptyset =A_0 \subset A_1 \subset A_2 \subset \cdots \subset A_k=E
\end{equation}
where $E$ is the ground set for the matroid $M$, let $\alpha(\flag)=(\alpha_1,\ldots,\alpha_k)$ 
be the composition of $n:=|E|$ defined by $\alpha_i:=|A_i|-|A_{i-1}|$.

\begin{proposition}
\label{ABS-computation}
The coefficient $c_\alpha$ in \eqref{M-expansion} 
is the number of flags $\flag$ of subsets of $E$ having $\alpha(\flag)=\alpha$
and for which each subquotient $\left( M|_{A_i} \right) / A_{i-1}$ splits completely.
\end{proposition}

The equivalence of these two interpretations of $c_\alpha$ is easy to understand.
Any $f: E \rightarrow \PP$ with $|f^{-1}(i)|=\alpha_i$ for $i=1,2,\ldots,k$
gives rise to a flag $\flag$ of subsets as in \eqref{typical-flag} with
$\alpha(\flag)=\alpha$, by letting $A_i:=f^{-1}(\{1,2,\ldots,i\})$.
In other words, $f$ takes on the constant value $i$ on each of the set differences
$A_i \setminus A_{i-1}$.  One can then readily see (e.g. from the greedy algorithm)
that $f$ will be $M$-generic if and only if each of the
subquotients $\left( M|_{A_i} \right) / A_{i-1}$ has only one base,
that is, if and only if each such subquotient splits completely.

A consequence of this equivalence is that \cite[Theorem 4.1]{AguiarBergeronSottile}
gives an alternate proof of Theorem~\ref{morphism-theorem} above.

Given this discussion, the existence of canonical {\it odd and even subalgebras} inside any combinatorial
Hopf algebra (see \cite[\S 5]{AguiarBergeronSottile}) naturally suggests the following question.

\begin{question}
\label{Eulerian-question}
What is the odd subalgebra of the combinatorial
Hopf algebra $\Mat$?  Does it contain any elements of the form $[M]$ for a single matroid $M$, or does it 
contain only nontrivial sums $\sum c_M [M]$? 
\end{question}

\noindent
For any such $[M]$ in the odd subalgebra of $\Mat$, it will follow from
\cite[Propositions 5.8e and 6.5]{AguiarBergeronSottile} that $F(M)$ will
lie in the {\it peak subalgebra} of $\Qsym$ (see \cite{AguiarBergeronSottile}
for definitions).

\end{remark}

\begin{example} \rm \ 
One can use Proposition~\ref{ABS-computation} to compute some more examples
of $F(M)$.  If $M$ is a rank $1$ matroid on ground set $E=\{1,2\}$ in which $1,2$ are
parallel elements, then there are exactly two flags $\flag$ having all subquotients
that split completely:
$$
\begin{aligned}
&\emptyset \subset \{1\} \subset E=\{1,2\} \\
&\emptyset \subset \{2\} \subset E=\{1,2\} \\
\end{aligned}
$$ 
Both of these flags have $\alpha(\flag)=(1,1)$ and hence $F(M)=2M_{1,1}$.

Similarly, if $M$ is a rank $1$ matroid on ground set $E=\{1,2,3\}$ in which
$1,2,3$ are all parallel, then there are two kinds of flags $\flag$ having all subquotients
that split completely: 
\begin{enumerate}
\item[$\bullet$]
$6=3!$ flags of the form
$$
\emptyset \subset \{a\} \subset \{a,b\} \subset E=\{a,b,c\}
$$ 
where $(a,b,c)$ is some permutation of $(1,2,3)$, all having $\alpha(\flag)=(1,1,1)$, and
\item
[$\bullet$]
$3$ flags of the form
$$
\emptyset \subset \{a\} \subset E=\{a,b,c\}
$$ 
where $a \in \{1,2,3\}$, all having $\alpha(\flag)=(1,2)$.
\end{enumerate}
Consequently, $F(M)=3M_{1,2} + 6M_{1,1,1}$.
\end{example}

\section{Behavior under matroid duality}
\label{matroid-duality-section}

Recall that if $M$ is a matroid on ground set $E$ with bases $\bases(M)$, then
its {\it dual} (or {\it orthogonal}) matroid $M^\ast$ has the same ground set $E$
but bases 
$$
\bases(M^\ast)=\{B^\ast: B \in \bases(M)\}
$$
where $B^\ast:=E\backslash B$ is called the 
{\it cobase} of $M^\ast$ corresponding to the base $B$ of $M$.

\begin{proposition}
\label{duality-proposition}
$$
F(M)  =\sum_{\alpha} c_\alpha M_\alpha
$$
if and only if
$$
F(M^\ast)=\sum_{\alpha} c_\alpha M_{\alpha^*}
$$
where $\alpha^* :=(\alpha_k,\alpha_{k-1},\ldots,\alpha_2,\alpha_1)$
is the reverse composition to $\alpha$.
\end{proposition}
\begin{proof}
We check that for any composition $\alpha \in \PP^k$,
the coefficient of $M_\alpha$ in $F(M)$ is
the same as the coefficient of $M_{\alpha^*}$ in $F(M^\ast)$.

The former coefficient counts the set of $M$-generic $f:E \rightarrow \PP$
for which $\xx_f=\xx^\alpha$.  The latter coefficient counts the
set of $M^\ast$-generic $f^\ast:E \rightarrow \PP$ for which $\xx_{f^\ast}=\xx^{\alpha^*}$.

We exhibit a bijection between these sets as follows.
If $B$ is a base of $M$ with cobase $B^\ast$ of $M^\ast$, then
the equation
$$
f(B) + f(B^\ast) = \sum_{e \in E} f(e)
$$ 
shows that $B$ is $f$-minimizing if and only if $B^\ast$ is $f$-maximizing.
Now define $f^\ast(e) := k+1-f(e)$, so that one has 
$$
f(B^\ast) + f^\ast(B^\ast) = (k+1) \; |B^\ast|  = (k+1)\left(\; |E|-r(M)\; \right).
$$
This equation shows that $B^\ast$ is $f$-maximizing if and only if $B^\ast$ is $f^\ast$-minimizing.
Since $\xx_f = \xx^\alpha$ if and only if $\xx_{f^\ast} = \xx^{\alpha^*}$,
the map $f \mapsto f^\ast$ restricts to the desired bijection.
\end{proof}

\section{$P$-partition expansion}
\label{P-partition-section}

Quasisymmetric functions were originally introduced by Gessel \cite{Gessel}
(building on work of Stanley) as enumerators for $P$-partitions.  
We review this here, and explain
how it leads to an expansion of $F(M)$ as a sum of $P$-partition
enumerators.

A {\it labelled poset} $(P,\gamma)$ on $n$ elements is a poset $P$ together with
a bijective labelling 
function $\gamma: P \rightarrow [n]:=\{1,2,\ldots,n\}$.  

A {\it $(P,\gamma)$-partition}
is a function $f: P \rightarrow \PP$ such that
$$
\begin{aligned}
f(p) \leq f(p') & \textrm{ if } p \leq p'  \\
f(p) <  f(p')    & \textrm{ if }p \leq p'\textrm{ and }\gamma(p) > \gamma(p')
\end{aligned}
$$
It will sometimes be more convenient for us to refer only
to a labelled poset $P$ on $[n]$ (suppressing the extra labeling function
$\gamma$), by which we mean a partial order $<_P$ on the set $[n]$.
Using this terminology, a $P$-partition is a function $f: [n] \rightarrow \PP$ satisfying
$$
\begin{aligned}
f(i) \leq f(i') & \textrm{ if }i \leq_P i' \\
f(i) < f(i')    & \textrm{ if }i \leq_P i' \textrm{ and }i >_{\ZZ} i'.
\end{aligned}
$$
For example, every
permutation $w=w_1 \cdots w_n$ of $[n]$ can be regarded as a labelled
poset on $[n]$ which is totally ordered:  $w_1 <_w \cdots <_w w_n$.

Let $\AAA(P,\gamma)$ denote the set of $(P,\gamma)$-partitions,
and let $F(P,\gamma,\xx):=\sum_{f} \xx_f$ be their weight enumerator:
$$
F(P,\gamma,\xx):=\sum_{f \in \AAA(P,\gamma)} \xx_f.
$$
A basic result of Stanley tells how $F(P,\gamma,\xx)$ expands
in terms of another basis for $\Qsym$ indexed by compositions $\alpha$, known as the
{\it fundamental quasisymmetric functions}
\begin{equation}
\label{L-into-M-expansion}
L_\alpha:=\sum_{\beta: \beta \textrm{ refines }\alpha} M_\beta.
\end{equation}

Say that a permutation $w=w_1 \ldots w_n$ in the symmetric group $S_n$
is a {\it linear extension} of $(P,\gamma)$ if $p < p'$ in $P$ 
implies $w^{-1}(\gamma(p)) < w^{-1}(\gamma(p'))$.  The {\it Jordan-H\"older
set} of $(P,\gamma)$ is the set $\LL(P,\gamma)$ of all linear extensions
of $(P,\gamma)$.  The {\it descent composition} for the permutation $w$
is the composition $\alpha(w)$ of $n$ which gives the lengths of the maximal
increasing consecutive subsequences ({\it runs}) of $w$.
It is not hard to check that, regarding $w$ as a totally ordered labelled poset
on $[n]$ as above, one has $ F( w, \xx)= L_{\alpha(w)}$.  The basic result about
$P$-partitions is the following expansion.

\begin{proposition} \cite[\S 4.5]{Stanley-EC1}, \cite[\S 7.19]{Stanley-EC2}, \cite[eqn. (1)]{Gessel}
\label{Stanley's-P-partition-result}
$$
\begin{aligned}
F(P,\gamma,\xx) &=\sum_{w \in \LL(P,\gamma)} F(w,\xx)\\
                &=\sum_{w \in \LL(P,\gamma)} L_{\alpha(w)}
\end{aligned}
$$
\end{proposition}

It turns out that every base $B$ of a matroid $M$ leads to a certain
labelled poset $P_B$, whose $P$-partition enumerator is relevant for expanding
$F(M)$;  see Theorem~\ref{poset-expansion} below.

Given a base $B$ of a matroid $M$ on ground set $E$, 
let $B^\ast=E\backslash B$ be the corresponding cobase
of $M^\ast$.  

For each $e \in B$ the {\it basic bond for $e$ in $B^\ast$} is the
set of $e' \in E$ for which $(B\backslash \{e\}) \cup \{e'\}$
is another base of $M$.  Dually, for each $e \in E-B (=B^\ast)$ the
{\it basic circuit for $e$ in $B$} is the
set of $e' \in E$ for which $(B\cup \{e\}) \backslash \{e'\}$
is another base of $M$.  By definition then, one has a symmetric relationship:
$e'$ lies in the basic bond for $e$ in $B^\ast$ if and only if
$e$ lies in the basic circuit for $e'$ in $B$.  Thus these relations
can be encoded by a bipartite graph with vertex set $E$, bipartitioned
as $E=B \sqcup B^\ast$.  
Define the poset $P_B$ to be the one whose
Hasse diagram is this bipartite graph, with edges directed upward
from $B$ to $B^\ast$.  

Say that a labelling $\gamma$ of a poset $P$
is {\it natural} (resp. {\it strict} or {\it anti-natural}) if
$\gamma(p) < \gamma(p')$ (resp. $\gamma(p) > \gamma(p')$)
whenever $p < p'$ in $P$.

\begin{theorem} 
\label{poset-expansion}
For any matroid $M$,
$$
F(M,\xx) = \sum_{B \in \bases(M)} F(P_B,\gamma_B,\xx)
$$
where $\gamma_B$ is any strict labelling of $P_B$.
\end{theorem}

\begin{proof}
We will show that $B$ is the unique $f$-minimizing base of $M$ for
some $f:E \rightarrow \PP$ if and only if $f$ lies in $\AAA(P_B,\gamma_B)$.

First assume that $f$ does not lie in $\AAA(P_B,\gamma_B)$, that is, there
exists some $e < e'$ in $P_B$ for which $f(e) \geq f(e')$.  By definition
of $P_B$, this means that $e$ lies in $B$, $e'$ does not lie in $B$, and
$B':=B\backslash \{e\} \cup \{e'\}$ is another base of $M$.  However
$f(e') \leq f(e)$ implies $f(B') \leq f(B)$, 
so that $B$ cannot be the unique $f$-minimizing base.

Now assume that $B$ is not the unique $f$-minimizing base of $M$.  This
means that there exists another base $B'$ of $M$ having $f(B') \leq f(B)$.
By convexity, we may assume that the pair $\{B,B'\}$ corresponds to an
edge of the {\it matroid base polytope} $Q(M)$, 
which is defined to be the convex
hull in $\RR^E$ of all characteristic $\{0,1\}$-vectors of bases of $M$
(see Section~\ref{decomposition-section} below).
A well-known fact from matroid theory \cite[\S2.2, Theorem 1]{GelSerg}
says that all edges of $Q(M)$ take the form $\{B,B'\}$ in which $B, B'$ differ by a single {\it basis exchange}:
there exists some $e \in B$ and $e' \in B'$ such that $B'=B\backslash \{e\} \cup \{e'\}$.
Thus $e < e'$ in $P_B$.  Since $f(B') \leq f(B)$ forces $f(e') \leq f(e)$, this means
$f$ is not in $\AAA(P_B,\gamma_B)$.
\end{proof}

\begin{remark}\rm \ 
Aguiar has pointed out that Theorem~\ref{poset-expansion} shows
the Hopf morphism $F: \Mat \rightarrow \Qsym$ factoring through the {\it Hopf algebra
$\mathcal P$ of (labelled) posets}, 
which is described (for unlabelled posets) in \cite[Example 2.3]{AguiarBergeronSottile}.
More precisely, one has a Hopf morphism 
$$
\begin{matrix} 
\Mat &\longrightarrow &{\mathcal P}\\
[M]  &\mapsto     &\sum_{B \in \bases(M)}[(P_B,\gamma_B)]
\end{matrix}
$$
and the usual $(P,\gamma)$-partition enumerator Hopf morphism 
$$
\begin{matrix}
{\mathcal P} & \longrightarrow & \Qsym\\
[(P,\gamma)] & \mapsto     & F(P,\gamma,\xx).
\end{matrix}
$$  
Then $F: \Mat \rightarrow \Qsym$ is the composite
of these two morphisms.
\end{remark}

\begin{corollary}
\label{last-coefficient}
Let $F(M) = \sum_{\alpha} c^M_\alpha L_\alpha$.  Then
\begin{enumerate}
\item[(i)] the coefficients $c^M_\alpha$ are nonnegative,
\item[(ii)] their sum $\sum_{\alpha} c^M_\alpha$ is $n!$ where $n:=|E|$, and
\item[(iii)] the coefficient $c^M_{1,1,\ldots,1}$ of $L_{1,1\ldots,1}$ is
the number of bases of $M$.
\end{enumerate}
\end{corollary}
\begin{proof}
Everything will follow from Proposition~\ref{Stanley's-P-partition-result}
and Theorem~\ref{poset-expansion}.  
Assertion (i) is immediate.

Assertion (ii) follows because each of the $n!$ linear orderings $e_1,\ldots,e_n$ of
$E$ is a linear extension for exactly one of the posets $P_B$, namely the one indexed
by the unique $f$-minimizing base $B$ when $f(e_1) < \cdots < f(e_n)$.

Assertion (iii) follows because any strictly (anti-naturally) labelled poset 
$(P,\gamma)$ has the reversing permutation $w_0 = n \ldots 3 2 1$
in $\LL(P,\gamma)$, and $w_0$ is the {\it only} permutation having
descent composition $(1,1,\ldots,1)$.
\end{proof}

Corollary~\ref{last-coefficient} gives a combinatorial interpretation for 
the coefficient $c^M_{1,1,\ldots,1}$. It would be
nice to have such an interpretation for every coefficient $c^M_\alpha$.  The next
result at least tells us how to interpret the coefficients ``at the
other end'' of the $L_\alpha$ expansion, namely $c^M_{\alpha}$ where $\alpha$ has at most two parts,
in terms of some basic matroid invariants of $M$.  
Recall that an element $e$ in $E$ is a {\it loop} in $M$ if
it appears in {\it no} bases of $M$, and it is a {\it coloop} (or {\it isthmus})
if it appears in {\it every} base of $M$.

\begin{proposition}
\label{first-coefficients}
Let $M$ be matroid having 
\begin{enumerate}
\item[$\bullet$] rank $r$, 
\item[$\bullet$] corank $r^\ast:=|E|-r$, 
\item[$\bullet$] number of loops equal to $\ell$,
\item[$\bullet$] number of coloops equal to $c$, and 
\item[$\bullet$] number of bases $b$.
\end{enumerate}

Then 
$$
F(M) = b\left( \sum_{j=0}^{\ell+c} \binom{\ell+c}{j} L_{(r+\ell-j,r^\ast-\ell+j)} \right)
                        + \sum_{\beta: \ell(\beta) \geq 3} c_\beta L_\beta
$$
for some nonnegative coefficients $c_\beta$.
Here $\ell(\beta)$ denotes the number of parts in the composition $\beta$.

Equivalently, if $\hat{M}$ is the matroid obtained from
$M$ by removing all loops and coloops, so that $\hat{M}$ has 
$$
\begin{aligned}
\text{ rank }\hat{r}&=r-c,\text{ and }\\
\text{ corank }\hat{r}^\ast&=r^\ast-\ell=|E|-r-\ell,
\end{aligned}
$$
then 
$$
\begin{aligned}
F(M) &= (L_1)^{\ell+c} F(\hat{M})  \\
     &= (L_1)^{\ell+c} \left( b \cdot  
         L_{(\hat{r},\hat{r}^\ast)} + \sum_{\gamma: \ell(\gamma) \geq 3} d_\gamma L_\gamma
       \right)
\end{aligned}
$$
for some nonnegative integer coefficients $d_\gamma$.

\end{proposition}
\begin{proof}
The second assertion follows from the $\ell+c=0$ case of the first,  applying the
multiplicative property $F(M_1 \oplus M_2) = F(M_1) F(M_2)$ to the decomposition of
$M$ as a direct sum of $\hat{M}$ with $\ell+c$ loops and isthmuses.

For the first assertion, we apply Theorem~\ref{poset-expansion}.  For each base $B$, the poset $P_B$ will
have height one, and decompose into three sets:  
\begin{enumerate}
\item[$\bullet$]
the set $A_1$ of $\ell+c$ loops and coloops, which are all both minimal and maximal in $P_B$, 
\item[$\bullet$]
the set $A_2$ of $r-c$ non-coloop elements in $B$, each of
which is minimal but not maximal in $P_B$, and 
\item[$\bullet$] 
the set $A_3$ of $r^\ast-\ell$ non-loop elements in $B^\ast$, each of which is
maximal but not minimal in $P_B$.
\end{enumerate}
We are free to choose the strict labelling $\gamma_B$ so that the elements in
$A_2$ all have the highest labels, the elements in $A_3$
all have the lowest labels, and the elements in $A_1$ have the labels
in between.  

How then can one choose a linear extension $w$ in $\LL(P_B,\gamma_B)$ so that
its descent composition $\alpha(w)$ has at most two parts?  This means that $w$ 
has at most two increasing runs, separated by a unique descent.  Because of our
chosen labelling of $B$, such a $w$ will have the
first run of length at least $r-c$, and the second run of length at least $r^\ast-\ell$.
Furthermore, for any integer $j$ in the range $[0,\ell+c]$, one can check that
there are $\binom{c + \ell}{j}$ ways to choose such a $w$ in $(P_B,\gamma_B)$ so that it
starts with an increasing
run of length $r-c+j$, followed by its unique descent, and then ends with
an increasing  run of length $r^\ast+c-j$:  one must place the elements of
$A_2$ together with any $j$ elements chosen from $A_1$ before the
unique descent, and place the elements of $A_3$ together with other $\ell+c-j$
elements of $A_1$ after the unique descent.  
\end{proof}

In particular, the previous proposition tells us the coefficient
of $L_\alpha$ in $F(M)$ when $\alpha$ has only $1$ part.
Recall from Section~\ref{Hopf-section} that a matroid $M$ is said to split completely if $M$ is
a direct sum of loops and isthmuses.

\begin{corollary}
\label{zeta-computation}
For $M$ a matroid on ground set $E$ of size $|E|=n$,
the expansion of $F(M)$ in the $L_\alpha$  (resp. $M_\alpha$) basis for $\Qsym$ has
the coefficient of $L_{(n)}$ (resp. $M_{(n)}$) equal to
$1$ if $M$ splits completely, and $0$ otherwise.
\end{corollary}
\begin{proof}
The assertion for the coefficient of $L_{(n)}$ follows from Proposition~\ref{first-coefficients}.
Then the assertion for the coefficient of $M_{(n)}$ follows from the expansion 
\eqref{L-into-M-expansion} of $L_\alpha$ into $M_\beta$'s.
\end{proof}

\section{Reciprocity and behavior under the antipode}
\label{antipodes}

  Part of the structure of a Hopf algebra is an involutive anti-automorphism known
as its {\it antipode}.  For the Hopf algebra of quasisymmetric functions,
the antipode $S:\Qsym \rightarrow \Qsym$ is known to be related to
combinatorial reciprocity results \cite{MalvenutoReutenauer, Stanley2}.
It turns out to have an interesting effect on $F(M)$, transforming it into a different sort of 
enumerator for weight functions $f: E \rightarrow \PP$.  We begin by reviewing how the
antipode relates to reciprocity.

The antipode $S: \Qsym \rightarrow \Qsym$ has the following effect on 
the $L_\alpha$-basis \cite[Corollary 2.3]{MalvenutoReutenauer}:
$$
S(L_\alpha) = (-1)^{|\alpha|} L_{\alpha^c}
$$ 
where $|\alpha|:=\alpha_1+\cdots + \alpha_k=n$ denotes the {\it weight} of the composition $\alpha$,
and $\alpha^c$ corresponds to the subset $T^c=[n-1] \backslash S$ if $\alpha$ corresponds to
the subset $T$ of $[n-1]$ ({\it i.e.}  $T$ is the set of partial sums of $\alpha$).  

Stanley's reciprocity theorem for $P$-partitions \cite[Theorem 4.5.7]{Stanley-EC1} 
tell us that if $\gamma, \bar\gamma$  are 
natural and strict labellings of the same poset $P$, then
\begin{equation}
\label{natural-strict-relation}
S(F(P,\gamma,\xx)) = (-1)^{|P|} F(P,\bar\gamma,\xx).
\end{equation}

Upon specializing $L_\alpha, L_{\alpha^c}$ to $\xx=1^m$, that is, 
$$
\begin{aligned}
x_1=\cdots=x_m=1,\\
x_{m+1}=x_{m+2} = \cdots =0,
\end{aligned}
$$
one obtains
$$
\begin{aligned}
L_\alpha(1^m) &= \binom{m-k+n}{n} \\
L_{\alpha^c}(1^m) &=\binom{m+k-1}{n}.
\end{aligned}
$$
where $\alpha=(\alpha_1,\ldots,\alpha_k)$.  Then the equality
$$
\binom{m-k+n}{n}=(-1)^n \binom{-m+k-1}{n}
$$
leads immediately to the following reciprocity fact (cf. \cite[\S 4]{Stanley2}).

\begin{proposition}
If two homogeneous quasisymmetric functions $F, \Fdual$ of degree $n$
are related by $S(F)=\Fdual$, then their specializations
$$
\begin{aligned}
\phi(m)&=F(1^m)\\
\fdual(m) &=\Fdual(1^m)
\end{aligned}
$$
satisfy
$$
\phi(-m) = \fdual(m).
$$
\end{proposition}

We can now identify the image of $F(M)$ under the antipode $S$ in $\Qsym$.

\begin{defn} \rm \
Define a power series in $x_1,x_2,\ldots$
$$
\Fdual(M,\xx):=\sum_{ f: E \rightarrow \PP}  |\{ f\textrm{-minimizing bases of }M\}| \; \xx_f .
$$
Also define a polynomial in $m$ 
$$
\begin{aligned}
\fdual(M,m)&:=\Fdual(M,1^m) \\
           &=\sum_{f: E \rightarrow [m]} |\{ f\textrm{-minimizing bases of }M\}|.
\end{aligned}
$$
\end{defn}

One could argue that these two enumerators $\Fdual(M,\xx), \fdual(M,m)$ 
are at least as natural to consider as our original $F(M,\xx), \phi(M,m)$.
For example, the expected number of $f$-minimizing bases of $M$ attained when using at most $m$ distinct
values for the weights is exactly $\frac{1}{m^n} \fdual(M,m)$.

\begin{theorem}
\label{reciprocity-theorem}
For any matroid $M$ on $n$ elements, 
$$
S( F(M,\xx) ) = (-1)^n \Fdual(M,\xx)$$
and consequently, 
$$
\phi(M,-m) = (-1)^n \fdual(M,m).
$$
\end{theorem}
\begin{proof}
Theorem~\ref{poset-expansion} implies
$$
\begin{aligned}
S( F(M,\xx) ) &= \sum_{B \in \bases(M)} S (F(P_B,\gamma_B,\xx) ) \\
              &= (-1)^n \sum_{B \in \bases(M)} F(P_B,\bar\gamma_B,\xx) \\
              &= (-1)^n \sum_{B \in \bases(M)} 
                 \sum_{\substack{f:E \rightarrow \PP \\ 
                 B\textrm{ is }f\textrm{-minimizing}}} \xx_f \\
              &= (-1)^n \sum_{f: E \rightarrow \PP} |\{ f\textrm{-minimizing bases of }M\}| \xx_f \\
              &= (-1)^n \Fdual(M,\xx).
\end{aligned}
$$
\end{proof}

Note that since $F(M,\xx), \Fdual(M,\xx)$ are related by the antipode $S$, they
carry equivalent information, a fact which is not completely obvious from their
definitions.  The same goes for $\phi(M,m)$ and $\fdual(M,m)$.

\section{Valuation property and application to polytope decompositions}
\label{decomposition-section}

  The goal of this section is to show that the 
matroid invariant $F(M)$ behaves like a valuation on
the associated matroid base polytopes $Q(M)$, 
and apply this to the subtle problem of detecting decompositions
of these polytopes.

By the {\it matroid base polytope} we mean the convex polytope
$$
Q(M):= \text{conv} \left\{\sum_{i\in B}e_{i} : B \text{~a base of~} M\right\},
$$
where $e_{i}$ denotes the $i^{th}$ standard basis vector in $\RR^{E}$.  This polytope
$Q(M)$ is a face of a polytope 
first studied by Edmonds \cite{Edmonds}, which took as vertices the indicator
functions of \emph{all} independent sets in $M$ (subsets of bases).
We are interested in the existence or non-existence of certain polytopal decompositions of $Q(M)$.

\begin{defn} \rm \
A {\it matroid base polytope decomposition} of
$Q(M)$ is a decomposition $Q(M) = \cup_{i=1}^t Q(M_{i})$ where 
\begin{enumerate}
\item[$\bullet$] each $Q(M_i)$ is a matroid base polytope for some matroid $M_i$,
and
\item[$\bullet$]
for each $i \neq j$, the intersection
$Q(M_{i}) \cap Q(M_{j})$ is a face of both $Q(M_i)$ and of $Q(M_j)$.
\end{enumerate}
\end{defn}

We call such a decomposition a {\it hyperplane split} of $Q(M)$ if $t=2$.
We say that $Q(M)$ is {\it decomposable} if it has a matroid base polytope
decomposition with $t \geq 2$, and {\it indecomposable} otherwise.
We say that the decomposition is {\it coherent} if the $Q(M_i)$ are exactly the maximal domains of 
linearity for some $\RR$-valued piecewise-linear convex
function on $Q(M)$.  For example, hyperplane splits are always coherent.

Coherent matroid base polytope decompositions arise in work of Lafforgue
\cite{Lafforgue1, Lafforgue2} on compactifications of the fine Schubert cell of the Grassmannian
corresponding to the matroid $M$, and in related work 
by Keel and Tevelev \cite[\S2.6]{KeelTevelev}, and by
Hacking, Keel and Tevelev \cite[\S 3.3]{HackingKeelTevelev}.
In particular, Lafforgue's work implies that
for a matroid $M$ represented by vectors in $\FF^r$, if $Q(M)$ is indecomposable, 
then $M$ will be {\it rigid}, that is, $M$ will have only finitely many
realizations, up to scaling and the action of $GL(r,\FF)$.

\subsection{Polar cones and valuations}

We will need a version of a theorem of Lawrence \cite[Theorem 16]{Law} (see also
 \cite[Corollary IV.1.6]{Barv}) about polarity, 
which can be proved by a minor adjustment
to the proof of \cite[Theorem IV.1.5]{Barv}.

Let $ \langle \cdot, \cdot \rangle$ denote the usual inner product on $\RR^n$.
If $A$ is a convex set in $\RR^n$, then denote by $[A]$ its indicator function and by
$I(A)$ the convex set
$$I(A):=  \{x \in \RR^n:  \langle x,y \rangle  > 0 \text{~for all~} y \in A \}.$$
Recall that a closed convex cone $K\subset \RR^{n}$ is said to be \emph{pointed} if it contains
no lines.  In this case, its \emph{polar cone}
$K^{\circ}:= \{x\in \RR^{n} : \langle x,y \rangle \le 0 \text{ for all } y\in K \}$
has a nonempty interior.  For a \emph{nonzero} pointed cone $K$, $I(K)$ is
the interior of $-K^{\circ}$.

We show that the function $A\mapsto I(A)$ acts as
a valuation on nonempty closed convex sets.

\begin{proposition}\label{polar_relation}
Let  $A_{1}, A_{2},\dots, A_{N}$ be a finite family of nonempty closed convex sets.
If
$$
\sum_i \alpha_i [A_i]=0
$$
for real numbers $\alpha_{1}, \alpha_{2},\dots$, then
$$
\sum_i \alpha_i [I(A_i)] =0.
$$
\end{proposition}

\begin{proof}
The proof is as in \cite{Barv}, except that in Theorem IV.1.5, one defines
$$F_{\epsilon}(x,y)=F(x,y)=
\begin{cases}
1 & \text{if~}  \langle x,y \rangle \le 0 \\
0 & \text{otherwise}.
\end{cases}
$$
In this case, the limiting argument of \cite{Barv} (and \cite{Law}) is not necessary.

As in \cite{Barv} the association
$${\mathcal{D}}: [A] \mapsto [I(A)]$$
is the specialization to indicator functions $[A]$ of a linear map
$${\mathcal{D}}: {\mathcal{C}}(\RR^d) \mapsto {\mathcal{C}}(\RR^d),$$
where ${\mathcal{C}}(\RR^d)$ is the algebra of indicator functions of closed convex sets in $\RR^d$ (see \cite[Defn. I.7.3]{Barv}). The map
${\mathcal{D}}$ may be defined as follows: for a function $g(x)$, the value of $({\mathcal{D}}g)(y)$ on a point $y \in \RR^d$ is given by
$$ \chi( g(x) ) - \chi( g(x)F(x,y) ).$$
Here $\chi$ denotes the Euler characteristic linear
functional on ${\mathcal{C}}(\RR^d)$;  its value on
a function $h(x) \in {\mathcal{C}}(\RR^d)$ is determined uniquely
from knowing that it takes the value $1$ on indicator functions of
closed convex sets.
\end{proof}

\subsection{Matroid polytopes and decompositions}

We now wish to apply this to a decomposition $Q(M)= \cup_i Q(M_{i}) $ of matroid base polytopes.

Necessarily the $M_{i}$ will be {\it weak images} (\emph{degenerations})
of $M$, that is, $\bases(M_{i})\subset \bases(M)$.
If $M_{1}=(E,\bases_{1})$ and $M_{2}=(E,\bases_{2})$ are matroids of the same rank on the
same set $E$, then we define $M_{1} \cap M_{2} := (E, \bases_{1}\cap \bases_{2})$.
We write $M_{1} \cap M_{2}= \emptyset$ if $ \bases_{1}\cap \bases_{2}=\emptyset$ (as opposed to
$\{\emptyset\}$).   Note that even when $M_{1} \cap M_{2}\not= \emptyset$,
it is not usually a matroid -- take, for example, the rank $2$ matroids
$M_1, M_2$ having bases 
$$
\begin{aligned}
\bases(M_1)&= \{13,14,23,24\},\\
\bases(M_2)&= \{12,13,23,24,34\}
\end{aligned}
$$
so that $\bases(M_{1}) \cap \bases(M_{2})=\{13,23,24\}$.  
However, when $Q(M_{1})$ and $Q(M_{2})$ meet along a 
common face (as in a matroid base polytope
decomposition), and that face is nonempty,
the intersection $M_1 \cap M_2$ {\it will} be a matroid.

\begin{proposition}
\label{matroid-intersection}
If $M_{1}$ and $M_{2}$ are matroids of the same rank $r$ and $Q(M_{1}) \cap Q(M_{2})$ is a
nonempty common face of $Q(M_{1})$ and $Q(M_{2})$, then $M_{1} \cap M_{2}$ is a matroid of
rank $r$ and
$$Q(M_{1}\cap M_{2})=Q(M_{1})\cap Q(M_{2}).$$
\end{proposition}
\begin{proof}
Nonempty faces of matroid base polytopes are matroid base polytopes
\cite[\S2.5 Theorem 2]{GelSerg}, and so the common face $Q(M_{1}) \cap Q(M_{2})$ must
be a matroid base polytope.  The vertices of $Q(M_{1}) \cap Q(M_{2})$
correspond to common bases of $M_{1}$ and $M_{2}$, that is, to elements of
$\bases_{1}\cap \bases_{2}$.
\end{proof}

Suppose $e_{B}=\sum_{i \in B} e_{i}$ is the vertex of $Q(M)$ corresponding to the base $B$ of $M$.
We denote by 
$K_{B}(M)$ the closed convex cone generated by the Minkowski sum (translate)
$Q(M) - \{e_{B}\}$.  Its polar
$K_{B}^{\circ}(M)$ is the {\it normal cone} to $Q(M)$ at $e_{B}$.

Notice that by the proof of Theorem \ref{poset-expansion}, the expansion of $F(M)$ given
there can be written
\begin{equation}\label{F-decomp}
\begin{aligned}
F(M,\xx) &= \sum_{B \in \bases(M)} F(K_B(M),\xx), \text{ where} \\
F(K_B(M),\xx) &= \sum_{f\in I\left(K_{B}(M)\right)}\xx_{f}
\end{aligned}
\end{equation}
With this, one can prove that $F(M)$ acts as a valuation over subdivisions of $Q(M)$.

\begin{theorem}\label{valuation}
The association $Q(M) \mapsto F(M)$ is a valuation on the class of matroid polytopes:  
if  $Q(M)$ can be subdivided into finitely many matroid polytopes 
$Q(M_{i})$, then
$$
F(M) = \sum_{j\ge1} (-1)^{j-1} \sum_{i_{1}<i_{2}<\cdots <i_{j}} 
F(M_{i_{1}}\cap M_{i_{2}}\cap \cdots \cap M_{i_{j}}),
$$
with the sum over $i_{1}<i_{2}<\cdots <i_{j}$ such that
$M_{i_{1}}\cap M_{i_{2}}\cap \cdots \cap M_{i_{j}}\ne \emptyset$.
\end{theorem}

\begin{proof}
Any decomposition of $Q(M)$ induces, for each $B\in \bases(M)$, a  decomposition
of $K_{B}(M)$ into $K_{B}(M_{i})$ where $B\in \bases(M_{i})$.  (For notational
convenience, we include all $B\in \bases$ and set
$K_{B}(M_{i})= \emptyset$ when $B\notin \bases(M_{i})$.)
This, in turn,
leads to an inclusion-exclusion relation (see, for example, \cite[Lemma I.7.2]{Barv})
\begin{align*}
\left[K_{B}(M)\right] &= \sum_{j}(-1)^{{j-1}}\sum_{i_{1}<i_{2}<\cdots<i_{j}}
				\left[K_{B}(M_{i_{1}})\cap\cdots\cap K_{B}(M_{i_{j}})\right]\\
				&=  \sum_{j}(-1)^{{j-1}}\sum_{i_{1}<i_{2}<\cdots<i_{j}}
				\left[K_{B}(M_{i_{1}}\cap\cdots\cap M_{i_{j}})\right],
\end{align*}
with the second equality following from Proposition \ref{matroid-intersection}.
Clearly, we can restrict these sums to those $i_{1}<i_{2}<\cdots<i_{j}$ for which
$B\in \bases(M_{i_{i}})\cap\cdots\cap\bases(M_{i_{j}})$, in which case 
$M_{i_{i}}\cap\cdots\cap M_{i_{j}}\ne \emptyset$.
Thus, by Proposition \ref{polar_relation}, we have the relation
\begin{align*}
\left[I\left(K_{B}(M)\right)\right] = \sum_{j}(-1)^{{j-1}}\sum_{i_{1}<i_{2}<\cdots<i_{j}}
				\left[I\left(K_{B}(M_{i_{1}}\cap\cdots\cap M_{i_{j}})\right)\right].
\end{align*}
The assertion now follows from \eqref{F-decomp}.
\end{proof}

It turns out that all of the terms with $j \geq 2$ in the summation of Theorem~\ref{valuation}
involve matroids which are {\it disconnected}.  This will allow us to deduce a corollary
(Corollary~\ref{sumresolution} below) which ignores these terms, and 
leaves a sum with {\it positive} coefficients.

To this end, recall that a nonempty subset $A \subseteq E$ is called a \emph{separator} of $M$ if it leads
to a direct sum decomposition of matroids:
$$
M = M|_{A} \oplus M|_{E\setminus A}
$$  
The whole ground set $E$ is itself a separator, 
and the collection of separators is closed under
intersection.  Hence $E$ can be written as a disjoint union of 
inclusion-minimal separators of $M$.  
Denote by $s(M)$ the number of minimal separators of $M$.  The following is
\cite[\S2.4, Proposition 4]{GelSerg}. 

\begin{proposition}
The dimension of the matroid polytope $Q(M)$ is $|E|-s(M)$.
\end{proposition}

Considering $\Qsym$ as a graded $\ZZ$-algebra,
its maximal (homogeneous) ideal is $\m = \oplus_{d \geq 1} \Qsym_d$.
Given an element $f \in \Qsym$, let $\overline{f}$ denote its image
in the quotient ring $\Qsym/\m^2$.

\begin{corollary}
If $E\ne \emptyset$ and the dimension of $Q(M)$ is less than $|E|-1$, then $F(M)$ lies in the
square $\m^{2}$ of the maximal ideal $\m$.  In other words,
$\overline{F(M)}=0$ in $\Qsym/\m^2$.
\end{corollary}
\begin{proof}
If $Q(M)$ has dimension less than $|E|-1$ then $s(M) > 1$, so there
exists at least one proper separator $A \subsetneq E$.
Since $F: \Mat \rightarrow \Qsym$ is 
an algebra morphism, one has $F(M) = F( M|_{A} ) F(M|_{E\setminus A})$, and hence
$F(M)$ lies in  $\m^{2}$.
\end{proof}

Since $Q(M_{1}\oplus M_{2}) = Q(M_{1}) \times Q(M_{2})$, to study decomposability of
matroid polytopes $Q(M)$, it is enough to restrict attention to \emph{connected} matroids $M$, that is,
those with $s(M)=1$.  For these, the maximal cells in any decomposition $Q(M)=\cup_i Q(M_i)$ will have dimension
$|E|-1$ and so will also correspond to connected matroids.  All their proper intersections, however,
will be lower-dimensional and so correspond to matroids with non-trivial separators.  

\begin{corollary}\label{sumresolution}
If a matroid polytope $Q(M)$ can be subdivided into finitely many matroid polytopes 
$Q(M_i)$, then in $\Qsym/\m^2$ one has 
$
\overline{F(M)} = \sum_{i} \overline{F(M_i)}.
$
\end{corollary}

This corollary interacts nicely with a result of Hazewinkel \cite[Theorem 8.1]{Haz}, 
confirming a conjecture of Ditters which says that the $\ZZ$-algebra structure on
$\Qsym$ is that of a {\it free} commutative algebra, that is,
a polynomial algebra.  Consequently, $\m/\m^2$ is a {\it free} (graded) $\ZZ$-module,
and hence each homogeneous component $(\Qsym/\m^2)_n$ is a $\ZZ$-lattice $\ZZ^{r_n}$ of some
finite rank\footnote{In fact, these ranks $r_n$ can be made more explicit in two ways.
First, they are determined uniquely by the power series relation
$$
\prod_{n \geq 1} \frac{1}{(1-t^n)^{r_n}} 
  = \mathrm{Hilb}(\Qsym,t)
  = 1 + t + 2t^2 + 4t^3 + \cdots = \frac{1-t}{1-2t}.
$$ 
Second, $r_n$ has a combinatorial interpretation explained in \cite[\S 4]{Haz},
as the number of words in the alphabet $\{1,2,\ldots\}$ of {\it total weight} $n$ which are
{\it star powers of elementary Lyndon words}.  

  In practice, we have done our computer calculations in $(\Qsym/\m^2)_n$
using $\{L_\alpha\}_{|\alpha|=n}$ as a $\ZZ$-basis for $\Qsym_n$,
and using $\{L_\beta L_\epsilon\}_{|\beta|+|\epsilon|=n}$ as a $\ZZ$-spanning set for $(\m^2)_n$.
To do this, one can expand $L_\beta L_\epsilon$ in terms of $L_\alpha$'s using 
Proposition~\ref{Stanley's-P-partition-result} above:
$L_\beta L_\epsilon = F(P,\gamma,\xx)$ for a labelled poset $(P,\gamma)$ which is the disjoint
union of two chains, one with descent composition $\beta$, the other with descent composition $\epsilon$.} $r_n$.

Thus for matroids $M$ of rank $r$ on ground set $E$ of size $n$,
to understand the potential matroid base polytope decompositions of $Q(M)$, it
helps to examine the {\it additive semigroup} structure generated by
the elements $\overline{F(M)}$ within the lattice $\ZZ^{r_n}$.

\begin{defn} \rm \
Say that $\overline{F(M)}$ is {\it decomposable} if there exist matroids $M_i$
with $\overline{F(M)} = \sum_{i=1}^t \overline{F(M_i)}$, and indecomposable
otherwise.  
Say that $\overline{F(M)}$ 
is {\it weak image decomposable} if it is decomposable with
each $M_i$ a weak image of $M$, that is,
$\bases(M_i) \subset \bases(M)$.
\end{defn}

In other words, $\overline{F(M)}$ is indecomposable if and only 
if it is an element of the unique {\it Hilbert basis} \cite[Chapter 13]{Sturmfels}
for the additive semigroup generated by the $\overline{F(M)}$.
Corollary~\ref{sumresolution} implies that $Q(M)$ is indecomposable unless
$\overline{F(M)}$ is weak image decomposable.  However, decomposability (and hence
also weak image decomposability) of
$\overline{F(M)}$ is easily checked using computer algebra packages that
can compute the {\it toric ideal} and/or the Hilbert basis for the additive
semigroup generated by the $\overline{F(M)}$ within $\Qsym/\m^2$;  
see \cite[Chapters 4 and 13]{Sturmfels}.

\begin{example} \rm \label{rank2}
A (loopless) rank $2$ matroid $M$ on $n$ elements is determined up to isomorphism by the partition
$\lambda(M)$ of $n$ that gives the sizes $\lambda_i$ of its parallelism classes.
Also, $M_1$ is a weak image of $M_2$, up to isomorphism,
if and only if the partition $\lambda(M_1)$ is refined by the partition $\lambda(M_2)$.  Note that
$M$ is connected if and only if $\lambda$ has at least $3$ parts.  Hence 
the connected weak images of $M$ correspond to all coarsenings
of $\lambda(M)$ with at least $3$ parts.

In particular, if $\lambda(M)$ has {\it exactly} $3$ parts then $\overline{F(M)}$
must be weak image indecomposable and $Q(M)$ must be indecomposable.  
In fact, by computer calculations, 
we have verified for $3 \le n \le 9$ that the rank $2$ matroids $M$ for which $\lambda(M)$ has
exactly $3$ parts form the Hilbert basis for the semigroup generated by the $\overline{F(M)}$, and 
those for which $\lambda(M)$ has more than $3$ parts all have $Q(M)$ decomposable 
(and hence $\overline{F(M)}$ weak image decomposable).  The following question was left open
in an earlier version of this paper, but has recently been resolved in the affirmative by
work of Luoto \cite[Corollary 6.7]{Luoto}:

\begin{question}
Fix $n$, and consider the semigroup generated by $\overline{F(M)}$ within $\Qsym/\m^2$
as one ranges over all matroids $M$ of rank $2$ on $n$ elements.
Is the  Hilbert basis for this semigroup indexed by
those $M$ for which $\lambda(M)$ has exactly $3$ parts?
\end{question}

We should mention that a convenient parametrization of all
rank $2$ matroid base polytope decompositions was given by Kapranov \cite[\S 1.3]{Kapranov},
who showed that all decompositions in this (rank $2$) setting 
can be achieved by a sequence of hyperplane splits.  
\end{example}

\begin{figure}
\begin{center}
\includegraphics[scale=0.58]{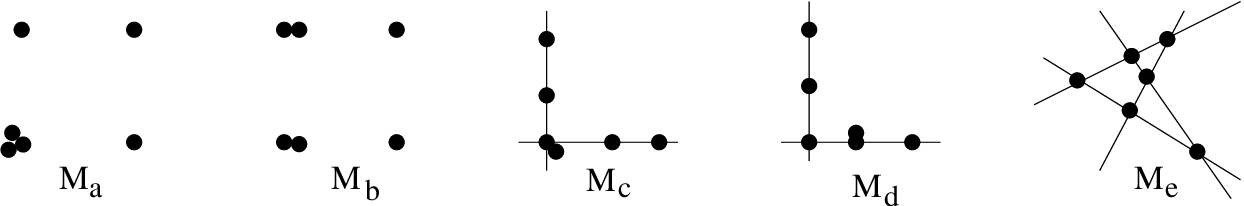}
\end{center}
 \caption{
\label{extremes-figure}
 The connected rank three matroids on 6 points for which $Q(M)$ is
indecomposable, represented
as affine point configurations in the plane.  The first
matroid contains a tripled point, and the second contains two pairs of doubled
points.
}
\end{figure}

\begin{example} \rm \label{rank3}
Considering all 15 connected rank 3 matroids $M$ with $n:=|E|=6$ (see, for example, \cite[Fig. 2]{GelSerg}), 
we found five for which $\overline{F(M)}$ is indecomposable.  
These are illustrated in Figure \ref{extremes-figure}.  

In particular, the two matroids $M_{1}$ and $M_{2}$ in
Figure~\ref{IsoTutte-figure}(a) satisfy $\overline{F(M_{1})}=\overline{F(M_{2})}$,
which can be written three different ways as sums of these indecomposables
\begin{align*}
\overline{F(M_{i})} &= \overline{F(M_{b})} + \overline{F(M_{c})} + 2 \overline{F(M_{d})}\\
&= 2 \overline{F(M_{a})} + \overline{F(M_{e})}\\
&= \overline{F(M_{a})} + 3 \overline{F(M_{d})} .  
\end{align*}
For $M_{2}$, all three of these additive decompositions
correspond to matroid base polytope decompositions of $Q(M_{2})$,
as does the first for $M_{1}$.
However, since $M_{a}$ is not a weak image of $M_{1}$, the second and third cannot
correspond to such decompositions of $Q(M_{1})$.
\end{example}

\begin{question}
Does $\overline{F(M)}$ being weak image decomposable in $\Qsym/\m^2$ imply that
$Q(M)$ is decomposable?
\end{question}

We see no reason, a priori, for this to hold,
but the matroids considered in Examples~\ref{rank2} and \ref{rank3} provide no
counterexamples. In fact, for all of the matroids $M$ in those examples,
one has $Q(M)$ indecomposable if and only if $\overline{F(M)}$ is indecomposable
if and only if $M$ is minimally connected 
({\it i.e.},  all weak images of $M$ have a nontrivial separator).

\begin{example} \rm \
It is worth noting that among the rank $3$ matroids with $n=6$ elements,
one finds the first matroid base polytope decompositions which are {\it not}
hyperplane splits.
For example, if $M$ is the rank 3 matroid on $E=\{1,2,3,4,5,6\}$ having every
triple
but $\{1,2,3\}$, $\{1,4,5\}$ and $\{3,5,6\}$ as bases, then both
$\overline{F(M)}$
and $Q(M)$ split into three indecomposable pieces, each isomorphic to the
matroid
$M_d$ in Figure \ref{extremes-figure}.
This subdivision of $Q(M)$ cannot be obtained via hyperplane splits.
\end{example}

\section{Comparison to the other matroid invariants}
\label{T-G-section}

One might ask how fine a matroid invariant is $F(M)$.  That is, how well does it
distinguish non-isomorphic matroids, say in comparison with well-studied
matroid invariants like the Tutte polynomial?

Certainly the kernel of the Hopf algebra map $F: \Mat \rightarrow \Qsym$
contains $p:=[M_\isthmus]-[M_\looop]$ by Example~\ref{tiny-examples}, 
and hence contains the smallest Hopf ideal $I$ generated
by $p$.  In fact, since $p$ is primitive (as it is of degree $1$),
the Hopf ideal $I$ which it generates coincides with the principal ideal
consisting of all multiples of $p$. Consequently $F$ factors through the quotient $\Mat /I$, that is,
through the Hopf algebra of matroids modulo ``loops = coloops''.

Beyond this inability to distinguish loops from coloops, one might
ask how discriminating $F(M)$ is. 
The next two examples show that it certainly doesn't distinguish all loopless and
coloopless matroids up 
to isomorphism (which would have been too much to ask), but it
at least does better than the well-known Tutte polynomial in some instances.

\begin{figure}
\begin{center}
\includegraphics[scale=0.3]{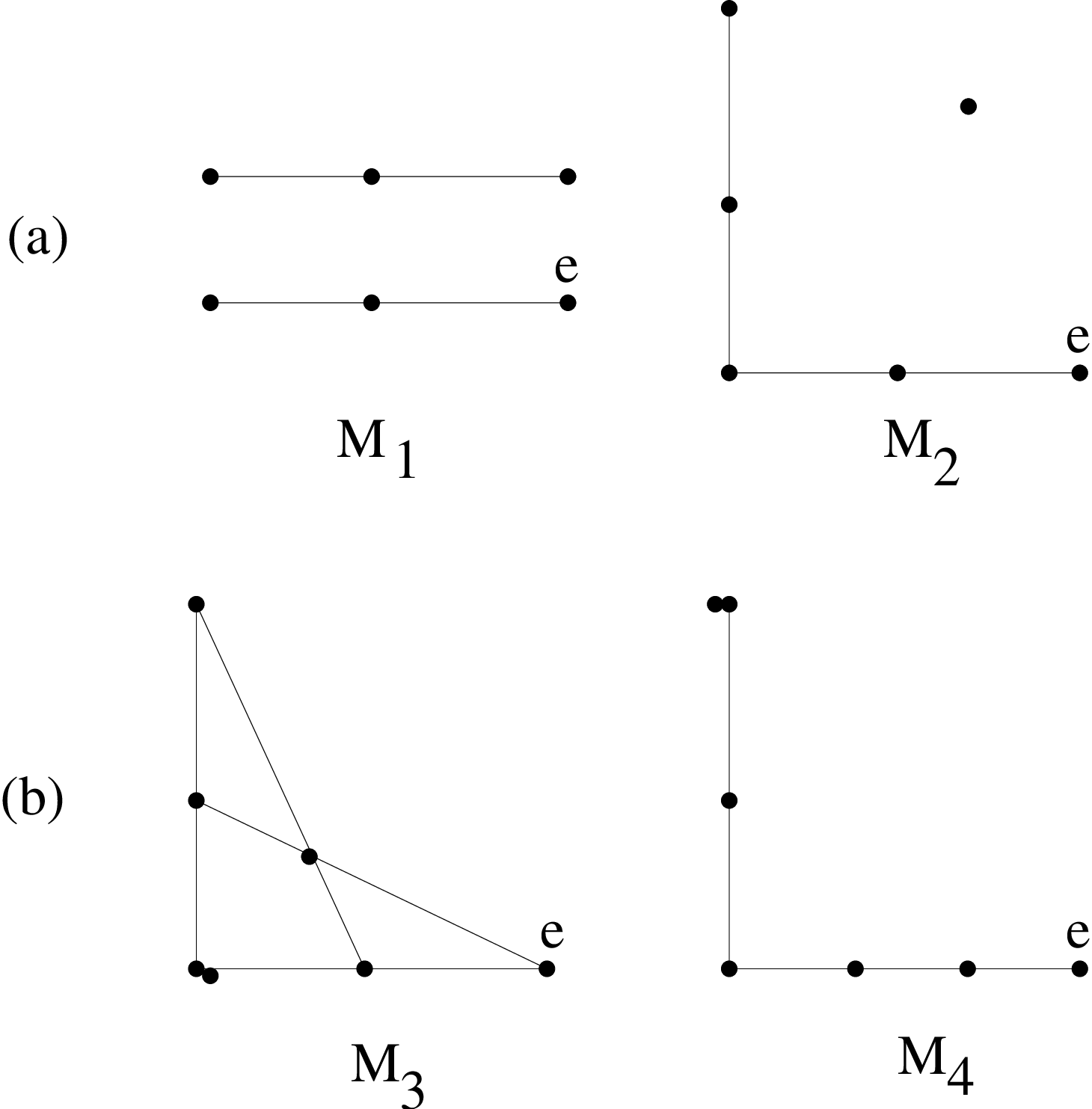}
\end{center}
 \caption{
\label{IsoTutte-figure}
  (a) Two matroids $M_1, M_2$, represented by affine point configurations
having the same Tutte polynomials and the same quasisymmetric functions.
(b) Two matroids $M_3, M_4$ having the same Tutte polynomials 
but different quasisymmetric functions.
}
\end{figure}

\begin{example} \rm \
\label{six-element-example}
Figure~\ref{IsoTutte-figure}(a) depicts two matroids, represented as affine point
configurations, having the same Tutte polynomial (because $M_1/e \cong M_2/e$ and 
$M_1 \backslash e \cong M_2 \backslash e$, where $e$ is the labelled point in
each case).  Direct computer calculation (using Theorem~\ref{poset-expansion}) shows that
$F(M_1)=F(M_2)$.
\end{example}

These two examples were borrowed from Brylawski and Oxley's survey on the Tutte polynomial
\cite[pp. 197]{White3}; they are the smallest examples of non-isomorphic matroids with
the same Tutte polynomial.

\begin{example} \rm \
\label{seven-element-example}
Figure~\ref{IsoTutte-figure}(b) depicts two matroids $M_3, M_4$ having the same Tutte polynomials 
(since $M_3/e \cong M_4/e$ and $M_3 \backslash e \cong M_4 \backslash e$), but
different quasisymmetric functions: it turns out that the coefficient of $L_{(1,3,3)}$
in $F(M_3)$ is $16$, while in $F(M_4)$ is $18$.  
\end{example}

These two examples were again taken from Brylawski and Oxley's survey 
\cite[pp. 133]{White3}, where they point out other features that $M_3, M_4$ share and do 
not share.

Note that Example~\ref{seven-element-example}  rules out the possibility of computing
$F(M)$ purely in terms of $F(M \backslash e), F(M/e)$.

\begin{question} \rm \  
Even though there is no direct deletion-contraction computation for $F(M)$ for
which one might have naively hoped, does this rule out other sorts of recursions?
\end{question}
In particular, one is tempted to try the following.
Theorem~\ref{poset-expansion} says
$$
\begin{aligned}
F(M,\xx) &= \sum_{B \in \bases(M)} F(P_B,\gamma_B,\xx) \\
         &= \sum_{B \in \bases(M): e \notin B} F(P_B,\gamma_B,\xx) +
            \sum_{B \in \bases(M): e \in B} F(P_B,\gamma_B,\xx).
\end{aligned}
$$
Can one better identify the two summands in this last equation?
Are they instances of some quasisymmetric functions that should be
associated to objects more general than matroids?

  Lastly, we mention an invariant $g_M(t)$  for a matroid $M$ (representable over $\QQ$) 
recently introduced by Speyer \cite{Speyer}, which shares some common features with
$F(M)$.  Among other of its properties, this invariant $g_M(t)$ is
\begin{enumerate}
\item[(i)] a polynomial in one variable t with integer coefficients 
(conjecturally nonnegative), 
\item[(ii)]
multiplicative under direct sums: $g_{M_1 \oplus M_2} = g_{M_1} g_{M_2}$,
\item[(iii)] invariant under duality of matroids: $g_M = g_{M^\ast}$,
\item[(iv)] additive under any decomposition of the matroid base polytope
$Q(M)=\cup_{i=1}^t Q(M_i)$, where $Q(M_1),\ldots,Q(M_t)$ are all of the {\it interior}
faces of the decomposition, in the sense that $g_M(t) = \sum_{i} g_{M_i}(t)$.
\end{enumerate}

\begin{question} \rm \ 
Is $g_M(t)$ related to (some specialization of) $F(M)$?  
\end{question}

%

\begin{remark} \rm \
In personal communication, Speyer has pointed out that all three invariants of
a matroid $M$ discussed in this section, behave either {\it valuatively} or 
{\it additively} under matroid base polytope decompositions:
\begin{enumerate}
\item[$\bullet$] the quasisymmetric function $F(M, \xx)$ behaves valuatively according to 
Theorem~\ref{valuation},
\item[$\bullet$] Speyer has checked that the Tutte polynomial $T_M(x,y)$ also behaves valuatively 
(via a small calculation using the {\it corank-nullity} formula for $T_M(x,y)$; 
see also \cite{ArdilaFinkRincon}), and 
\item[$\bullet$] his invariant $g_M(t)$ behaves additively by property (iv) above.
\end{enumerate}

\noindent
Speyer then used this to explain why all three invariants take the same value for the
two matroids $M_1, M_2$ shown in Figure~\ref{IsoTutte-figure}(a):  either of the matroid
base polytopes $Q(M_i)$ for $i=1,2$ can be obtained from the hypersimplex $Q(U(3,6))$
associated to the uniform matroid of rank $3$ on the same six elements, by splitting off
with hyperplanes (in any order) two other polytopes $Q(M_i'), Q(M_i'')$, that is,
$$
Q(U(3,6)) = Q(M_i) \cup Q(M_i') \cup Q(M_i'').
$$
Furthermore, the $M_i', M_i''$ are {\it all isomorphic} as matroids 
$$
M_1' \cong M_2' \cong M_1'' \cong M_2''.
$$
and have {\it isomorphic intersections}
$$
M_1 \cap M_1' \cong M_2 \cap M_2' \cong M_1 \cap M_1'' \cong M_2 \cap M_2''.
$$
As a consequence, a matroid invariant $f(M)$ will have
$$
\begin{aligned}
& f(U(3,6)) = \\
&
\begin{cases}
f(M_i) + f(M_i') + f(M_i'') - f(M_i \cap M_i') - f(M_i \cap M_i'') &\text{ if }f\text{ is valuative},\\
f(M_i) + f(M_i') + f(M_i'') + f(M_i \cap M_i') + f(M_i \cap M_i'') &\text{ if }f\text{ is additive}
\end{cases}
\end{aligned}
$$
for $i=1,2$.  In either case, this forces $f(M_1)=f(M_2)$.

This strongly suggests trying to define a ``universal'' valuative invariant of
matroids, following McMullen's {\it polytope algebra}, 
and in particular his section \cite[\S 20]{McMullen} dealing with valuations invariant under a 
finite group action.  Build an abelian group starting with the free abelian group
on basis elements $[M]$ indexed by matroids $M$, imposing the valuation relation
for each matroid base polytope decomposition of $Q(M)$, and the relation $[M]=[M']$ if $M$ and
$M'$ are isomorphic as matroids\footnote{As McMullen points out in \cite[\S 20]{McMullen}, 
imposing invariance under finite group action (such as matroid isomorphism)
seems to require sacrificing the multiplicative structure in the polytope algebra
coming from Minkowski addition.  It also appears that in our situation 
one must sacrifice translation-invariance, and the structure coming from dilatations,
as the vertices of each matroid base polytope $Q(M)$ 
are required to be $\{0,1\}$-vectors whose coordinates sum to the rank $r(M)$.}.
Valuative matroid invariants are exactly the linear functionals on this abelian group.
%

\begin{problem} \rm \ 
Study the structure of this abelian group.  Are there special classes of special matroids
which generate it?
\end{problem}

\noindent
For example, a conjecture of Speyer \cite[Conjecture 11.3]{Speyer} would follow
if this abelian group were generated by the classes $[M]$ where $M$ runs over
all direct sums of {\it series-parallel matroids}.

\end{remark}

\section{Generalization to generalized permutohedra}
\label{generalized-permutohedra-section}

It turns out that the proofs of Proposition~\ref{quasisymmetry}, Theorem~\ref{poset-expansion},
Corollary~\ref{last-coefficient}, Theorem~\ref{reciprocity-theorem}, and Theorem~\ref{valuation}
generalize in a straightforward way to give results about a general class of convex polytopes
studied recently by Postnikov \cite{Postnikov};  see also
\cite{SturmfelsEtAl} and \cite{PRW}.

Given a convex
polytope $Q$ in $\RR^n$, the following conditions are well-known to
be equivalent \cite[Proposition 7.12]{Ziegler}:
\begin{enumerate}
\item[$\bullet$] Every edge of $Q$ lies in one of the directions
$\{e_i - e_j: 1 \leq i \neq j \leq n\}$.
\item[$\bullet$] The normal fan of $Q$ in $(\RR^n)^*$ is refined by the
usual {\it braid arrangement} (or {\it type $A_{n-1}$ Weyl chamber fan}).
\item[$\bullet$] The polytope $Q$ is a Minkowski summand of some realization of the 
{\it permutohedron} as a Minkowski sum of line segments (possibly of different lengths)
in the directions 
$\{e_i - e_j : 1 \leq i < j \leq n\}$.
\end{enumerate}
Say that $Q$ is a {\it generalized $n$-permutohedron} when any of these
equivalent conditions hold\footnote{Actually, the definition of generalized permutohedra given
in \cite{Postnikov} looks slightly different, but is shown to be equivalent to these
conditions in \cite[Appendix]{PRW}.}.

\begin{example} \rm \ 
\label{matroid-base-polytope-example}
Given a matroid $M$ on ground set $E=[n]$, the matroid base polytope $Q(M)$ defined
in Section \ref{decomposition-section} is a generalized   
 $n$-permutohedron \cite[\S2.2, Theorem 1]{GelSerg}, a
fact that played a crucial role in the proof of Theorem~\ref{poset-expansion}.
\end{example}

Given a polytope $Q$ in $\RR^n$, say that a function $f:[n] \rightarrow \PP$
(which we think of as giving an element of $(\RR^n)^*$) is {\it $Q$-generic}
if $f$ maximizes over $Q$ uniquely at a vertex.  In other words, $f$ lies in
the {\it interior} of an $n$-dimensional cone in the normal fan for $Q$.
One can then prove the following:

\begin{theorem}
\label{permutohedron-summand-theorem}
If $Q$ is a generalized $n$-permutohedron in $\RR^n$, then
\begin{enumerate}
\item[(i)] the power series
$$
F(Q,\xx):=\sum_{\substack{Q\textrm{-generic} \\f:[n] \rightarrow \PP}} \xx_f
$$
is quasisymmetric, with 
\item[(ii)]
an expansion in terms of $P$-partitions enumerators as
$$
F(Q,\xx) = \sum_{ \textrm{vertices }v\textrm{ of }Q} F(P_v,\gamma_v, \xx)
$$
where $(P_v,\gamma_v)$ are certain strictly labelled posets indexed
by the vertices of $Q$.
\item[(iii)] 
Furthermore, the coefficients $c^Q_\alpha$ in its expansion
$F(Q,\xx) =\sum_{\alpha} c^Q_\alpha L_\alpha$  
  \begin{enumerate}
   \item[(a)] are nonnegative,
   \item[(b)] sum to $n!$, and
   \item[(c)] have $c^Q_{1,1,\ldots,1}$ equal to the number the number of vertices of $Q$.
  \end{enumerate}
\item[(iv)] The antipode $S$ on $\Qsym$ satisfies
$$
S( F(Q,\xx) ) = (-1)^n \Fdual(Q,\xx)
$$
where
$$
\Fdual(Q,\xx) := \sum_{f:[n] \rightarrow \PP} |\{ f\textrm{-minimizing vertices of }Q\}|.
$$
\item[(v)]  The two polynomials $\phi(Q,m), \fdual(Q,m)$ in the variable $m$ defined
by specializing $F(Q,\xx), \Fdual(Q,\xx)$ to $\xx=1^m$ satisfy
$$
\phi(Q,-m) = (-1)^n \fdual(Q,m).
$$
\item[(vi)] Suppose $Q = \cup_i Q_i$ is a decomposition of $Q$ into finitely many
  permutohedron summands $Q_i$
{\bf introducing no new vertices}\footnote{The authors are
  grateful to Mario Sanchez for pointing out (in Feb. 2020) the necessity of this assumption, which holds for all
  matroid base polytope decompositions, and is used in the proof of 
  Theorem~\ref{valuation}.  This assumption was missing in
  previous versions of this paper,
  including the journal version.}, that is, the vertices of
each $Q_i$ are a subset of the vertices of $Q$,
and in which $Q_i \cap Q_j$ is a common face of $Q_i$ and
$Q_j$ for all $i,j$.
Then 
$$
F(Q,\xx) = 
\sum_{j\ge1} (-1)^{j-1} \sum_{i_{1}<i_{2}<\cdots <i_{j}} 
F(Q_{i_{1}}\cap Q_{i_{2}}\cap \cdots \cap Q_{i_{j}}),
$$
where the sum is over those terms in which $Q_{i_{1}}\cap Q_{i_{2}}\cap \cdots \cap Q_{i_{j}}$
is nonempty.
\end{enumerate}
\end{theorem}

\noindent
In fact, the posets $P_v$ appearing in the theorem have a very simple description:
$P_v$ is the transitive closure of the binary relation on $[n]$ which
has $i <_{P_v} j$ if there exists an edge of $Q$ of the form
$\{v,v'\}$ with $v'-v = e_j - e_i$.

In the remainder of this section, we discuss three naturally occurring
families of generalized $n$-permutohedra that have occurred in the literature.

\begin{problem} \rm \  
\label{gen-perm-qsym}
Study the quasisymmetric functions $F(Q,\xx)$ associated with any
of these families of generalized permutohedra $Q$.
\end{problem}

\subsection{Graphic zonotopes and Stanley's chromatic symmetric function.} \ \\
Let $G$ be a simple graph on vertex set $[n]$.  Let $Z_G$ denote the Minkowski sum of
line segments in the directions 
$$
\{e_i -e_j: \{i,j\} \text{ is a an edge of }G\}.
$$ 
Then $Z_G$ is a generalized $n$-permutohedron;
the $n$-permutohedron itself equals $Z_{K_n}$ where $K_n$ is the complete graph
on $n$ vertices.  
It is easy to see that a function $f:[n] \rightarrow \PP$ is $Z_G$-generic
if and only if it is a {\it proper coloring} of the vertex set $[n]$ of $G$.
One concludes that $F(Z_G,\xx)$ is the same as the {\it chromatic symmetric
function} $X_G(x_1,x_2,\ldots)$ introduced by Stanley \cite{Stanley2}, and studied
further by others in recent years.  Many of the results of this paper were
inspired by his work, and in particular Theorem~\ref{permutohedron-summand-theorem}
generalizes a few of the facts about $X_G$.

  It is also known  (see \cite[Example 4.5]{AguiarBergeronSottile}) 
that the map $G \mapsto X_G$ can be interpreted as a Hopf morphism  
between a certain Hopf algebra of graphs and the 
Hopf algebra $\Lambda$ of {\it symmetric functions} inside the quasisymmetric functions
$\Qsym$.  As far as we know, this morphism is of a different nature than our Hopf morphism
$F: \Mat \rightarrow \Qsym$.

\subsection{Polymatroids and flag matroids.} \ \\
  Example~\ref{matroid-base-polytope-example} alludes to
a famous resut of Gelfand, Goresky, MacPherson, and Serganova, characterizing
matroids in terms of their matroid base polytopes, which we rephrase slightly here.

\begin{theorem} 
(see \cite[\S2.2, Theorem 1]{GelSerg}, \cite[Theorem 1.11.1]{BorovikGelfandWhite}) 

Let $\bases$ be a collection of $r$-subsets of $[n]$, and $Q$ the
convex hull of their characteristic vectors in $\{0,1\}^n \subset \RR^n$.
Then $\bases$ is the collection of bases $\bases(M)$ for some matroid $M$ on ground set $E=[n]$
(and $Q=Q(M)$ is the associated matroid base polytope) if and only if 
$Q$ is a generalized $n$-permutohedron.

\end{theorem}

This led Gelfand, Goresky, MacPherson, and Serganova to the notion of
{\it Coxeter matroids} \cite{BorovikGelfandWhite}.  A Coxeter matroid is
the result of taking the characterization in the previous theorem and
\begin{enumerate}
\item[$\bullet$]
replacing $r$-subsets of $[n]$, which can be thought of as 
the cosets of maximal parabolic subgroups in the Coxeter
group of type $A_{n-1}$, with cosets of an arbitrary parabolic subgroup in an
arbitrary finite Coxeter group,
\item[$\bullet$] replacing the characteristic vectors of $r$-subsets 
with $W$-translates of sums of fundamental dominant weights,
\item[$\bullet$] replacing generalized $n$-permutohedra with Minkowski summands of
the zonotopes generated by other root systems.
\end{enumerate}

When the Coxeter group is of type $A_{n-1}$, considering arbitrary parabolic subgroups
instead of just maximal ones leads to the notion of a {\it flag matroid}, and its
{\it flag matroid base polytope}.  These will be generalized $n$-permutohedra
generalizing
the matroid base polytopes, whose vertices are vectors in $\NN^n$ that no longer necessarily sum to $r$,
but obey certain constraints on the sizes of their coordinates;  see \cite[\S 1.11]{BorovikGelfandWhite}.

Generalizing in another direction, 
a {\it discrete polymatroid base polytope of rank }$r$ (see \cite{HerzogHibi})
is a generalized $n$-permutohedron, 
each of whose vertices has nonnegative integer coordinates summing to $r$.
These polytopes were introduced by Edmonds \cite{Edmonds} in the context of combinatorial optimization.

\subsection{Graph-associahedra.} \ \\
Building on work of others (De Concini-Procesi,  Davis-Januszkiewicz-Scott, and
Carr-Devadoss),
Postnikov \cite{Postnikov} showed that the generalized $n$-permutohedra contain
an interesting subclass of polytopes
called {\it graph-associahedra}, indexed by simple graphs $G$ on vertex set $[n]$.  Within
this subclass, the associahedra and cyclohedra correspond to the cases where the graphs $G$ are
paths and cycles, respectively.

\section{Appendix: surjectivity and new bases for $\Qsym$}
\label{appendix}
\subsection{Sketch of surjectivity} 
\label{surjectivity-section}
The goal of this appendix is to prove the following.

\begin{theorem}
\label{surjectivity-theorem} 
The Hopf algebra morphism $F: \Mat \rightarrow \Qsym$ is surjective when one extends the scalars
to a field $\FF$ of characteristic zero.
\end{theorem}

\noindent
We observe here that the morphism $F$ is definitely {\it not} surjective without extending
scalars.  The image of the map $\Mat_2 \overset{F}{\rightarrow} \Qsym_2$ on homogeneous components
of degree $2$ is a sublattice of index $2$ within
$\Qsym_2$: there are only four non-isomorphic matroids on $2$ elements,
whose images under $F$ are all either of the form $L_{1,1}+L_{2}$ or $2L_{1,1}$.

Our approach will be to define, for each degree $n$,
a family of $2^{n-1}$ matroids on ground set $E=[n]$, whose images under $F$ span
$QSym_n$ with rational coefficients. It turns out that it will suffice to take
a subfamily of a family of $2^n$ matroids 
which were called {\it freedom matroids} in \cite{CrapoSchmitt1}, and
which we will call {\it $PI$-matroids} here.  They were considered in the context
of face enumeration in \cite{Liu} and in \cite{BER}, where they arose
in the context of combinatorial operators on zonotopes.

Given a matroid $M$,
let $I(M):=M \oplus M_\isthmus$ be a single-element extension of $M$ by
an isthmus.  Let $P(M)$ be a single-element extension of $M$ which is the
{\it principal extension of $M$ along the improper flat}, that is, one adjoins
a new element $e$ to the ground set, which is generic while obeying the constraint
that it does not increase the rank.

Say that $M$ is a {\it $PI$-matroid} if it can be
obtained from the empty matroid $M_\varnothing$ on 
$E=\varnothing$ by performing a sequence of repeated $M \mapsto I(M)$
and/or $M \mapsto P(M)$ operations.  It happens that every matroid with $|E|\leq 3$
is isomorphic to a $PI$-matroid.

Let $0\{0,1\}^{n-1}$ denote the collection of all 
binary strings $\sigma \in \{0,1\}^n$ that begin with a $0$.
Given $\sigma$ in $0\{0,1\}^{n-1}$, let $M_\sigma$ be the $PI$-matroid 
built from this sequence beginning with an empty matroid, where one performs 
the $I$ operation for each
$0$ and the $P$ operation for each $1$ in $\sigma$.  For example, the sequence $01111$
would build the $PI$-matroid $M_{01111}$ of rank $1$ consisting of $5$ parallel elements).

We will prove the following refinement of Theorem~\ref{surjectivity-theorem}.

\begin{theorem} 
\label{refined-surjectivity-thm}
The quasisymmetric functions 
$$
\{ F(M_\sigma) : \sigma \in 0\{0,1\}^{n-1} \}
$$
span $\Qsym_n \otimes \FF$ whenever $n!$ is invertible in $\FF$.
\end{theorem}

\begin{remark} \rm \  
The operation $M \mapsto I(M)$ which adds an isthmus to $M$ 
has a predictable effect on $F(M)$: 
$$
F( I(M)) = L_1 \cdot F(M).
$$
Seeing this, one might hope to approach Theorem~\ref{refined-surjectivity-thm} 
by understanding how $F( P(M))$ relates
to $F(M)$.  Unfortunately, $F( P(M))$ does not depend solely on $F(M)$ via
some operation in $\Qsym$.  For example, the two matroids
$$
\begin{aligned}
M_1 & := M_\isthmus \oplus M_\isthmus, \textrm{ and }\\
M_2 & := M_\isthmus \oplus M_\looop
\end{aligned}
$$
have $F(M_1)=F(M_2) (= L_1^2)$, however  
$$
\begin{aligned}
F( P(M_1)) &= 3L_{2,1}+3L_{1,1,1}, \textrm{ while}\\
F( P(M_2) ) &= 2L_{2,1} + 2L_{1,2} + 2L_{1,1,1}.
\end{aligned}
$$
\end{remark}

Instead, the proof of Theorem~\ref{refined-surjectivity-thm} (and hence
Theorem~\ref{surjectivity-theorem}) proceeds in three steps,
carried out over this and the next two subsections.

\medskip

\noindent \emph{Step 1.}  Introduce a family of posets $R_\sigma$ on $[n]$, also indexed
by $0\{0,1\}^{n-1}$, and show that the expansion of the
$F(M_\sigma,\xx)$ in terms of the {\it strictly labelled} $P$-partition enumerators for the $R_\sigma$
is triangular in some ordering.  Furthermore, the diagonal coefficients in this expansion
are products of binomial coefficients that all divide $n!$.
\medskip

\noindent \emph{Step 2.}  Introduce another family of labelled posets $Q_\sigma$ on $[n]$,
also indexed by $0\{0,1\}^{n-1}$, which are easily seen to form a $\ZZ$-basis for $\Qsym$,
and have some nice properties.
\medskip

\noindent \emph{Step 3.}   Show that the expansion of the {\it naturally labelled}
$P$-partition enumerators for the $R_\sigma$ in terms of the $P$-partition
enumerators for the $Q_\sigma$ is unitriangular with respect to some ordering.
From Step 2 it then follows that the former $P$-partition enumerators also give a $\ZZ$-basis
for $\Qsym$.  

\medskip

\noindent 
Then by equation \eqref{natural-strict-relation}, 
the strict $P$-partition enumerators of the $R_\sigma$ also give a $\ZZ$-basis for
$\Qsym$, and together with Step 1 this proves  Theorem~\ref{refined-surjectivity-thm}.

\medskip

Step 1 is completed in the remainder of this subsection, 
while Steps 2 and 3 are achieved in Subsections \ref{Q-basis-section} and \ref{R-Q-expansion-section}.  As mentioned in an earlier footnote, Luoto \cite[\S 7.4]{Luoto} has recently
found an alternative to Steps 2 and 3, by expanding the $R_\sigma$ basis elements
unitriangularly in terms of his ``matroid-friendly'' basis for $\Qsym$.

\medskip

Given $\sigma$ in $0\{0,1\}^{n-1}$, let $R_\sigma$ be the
labelled poset of height 1 (or 0) on $[n]$ 
having $i <_{Q_\sigma} j$ if $\sigma_i = 0, \sigma_j = 1$ and $i<j$.

Each such $\sigma$ also defines a partition of the set $[n]$ into
intervals that we will call the {\it blocks} $A_1,\ldots,A_t$ of $\sigma$, 
by breaking $[n]$ between the
positions $i,i+1$ where $(\sigma_i,\sigma_{i+1})= (1,0)$.
We also define a vector $(z_1,\ldots,z_t)$ associated to $\sigma$ as follows: 
$z_i$ is the number of positions $j$ in
the block $A_i$ for which $\sigma_j=0$.  It is not hard to see that
one can recover $\sigma$ uniquely from the blocks $(A_1,...,A_t)$ and
the values $(z_1,\ldots,z_t)$.

\begin{example} \rm \ 
Let $n=7$ and let $\sigma$ be the string in $0\{0,1\}^{n-1}$ given by
$$
\begin{matrix}
\sigma = &0&1&0&0&1&1&1&0 \\
         &1&2&3&4&5&6&7&8.
\end{matrix}
$$
Then $R_\sigma$ is the labelled poset on $[8]$ in which the
minimal elements are $1,3,4,8$, the maximal elements
are $2,5,6,7$ (and $8$), and the order relations are 
$$
\begin{aligned}
1 &<2,5,6,7 \\
3,4 &<5,6,7\\
\end{aligned}
$$
as illustrated in Figure~\ref{R-poset-figure}.

\begin{figure}
\begin{center}
\includegraphics[scale=0.30]{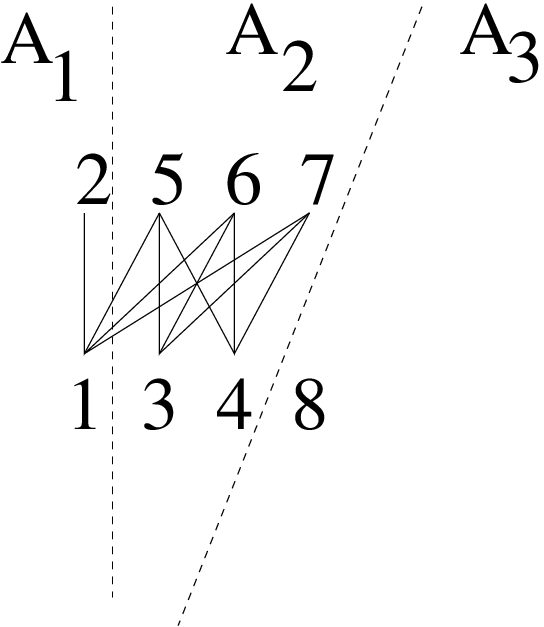}
\end{center}
 \caption{
\label{R-poset-figure}
The poset $R_{01001110}$, along with its associated blocks $(A_1,A_2,A_3)$.
}
\end{figure}

Also, $\sigma$ has associated to it the blocks $(A_1,A_2,A_3)=(12,34567,8)$,
and vector $(z_1,z_2,z_3)=(1,2,1)$.
The blocks $A_i$ are separated by dotted lines in Figure~\ref{R-poset-figure}.
\end{example}

It should be clear that the posets 
$R_\sigma$ are characterized up to isomorphism by the following {\it stable/shifted labelling} property.
\begin{proposition}
\label{shifted-characterization}
A labelled poset $P$ on $[n]$ is isomorphic to $R_\sigma$ for some $\sigma$ if and only if it has height at most
one, and can be relabelled so that each minimal (resp. maximal) element has its upward (resp. downward)
neighbors in $P$ forming a final (resp. initial) segment of $[n]$.
\end{proposition}

\begin{figure}
\begin{center}
\includegraphics[scale=0.40]{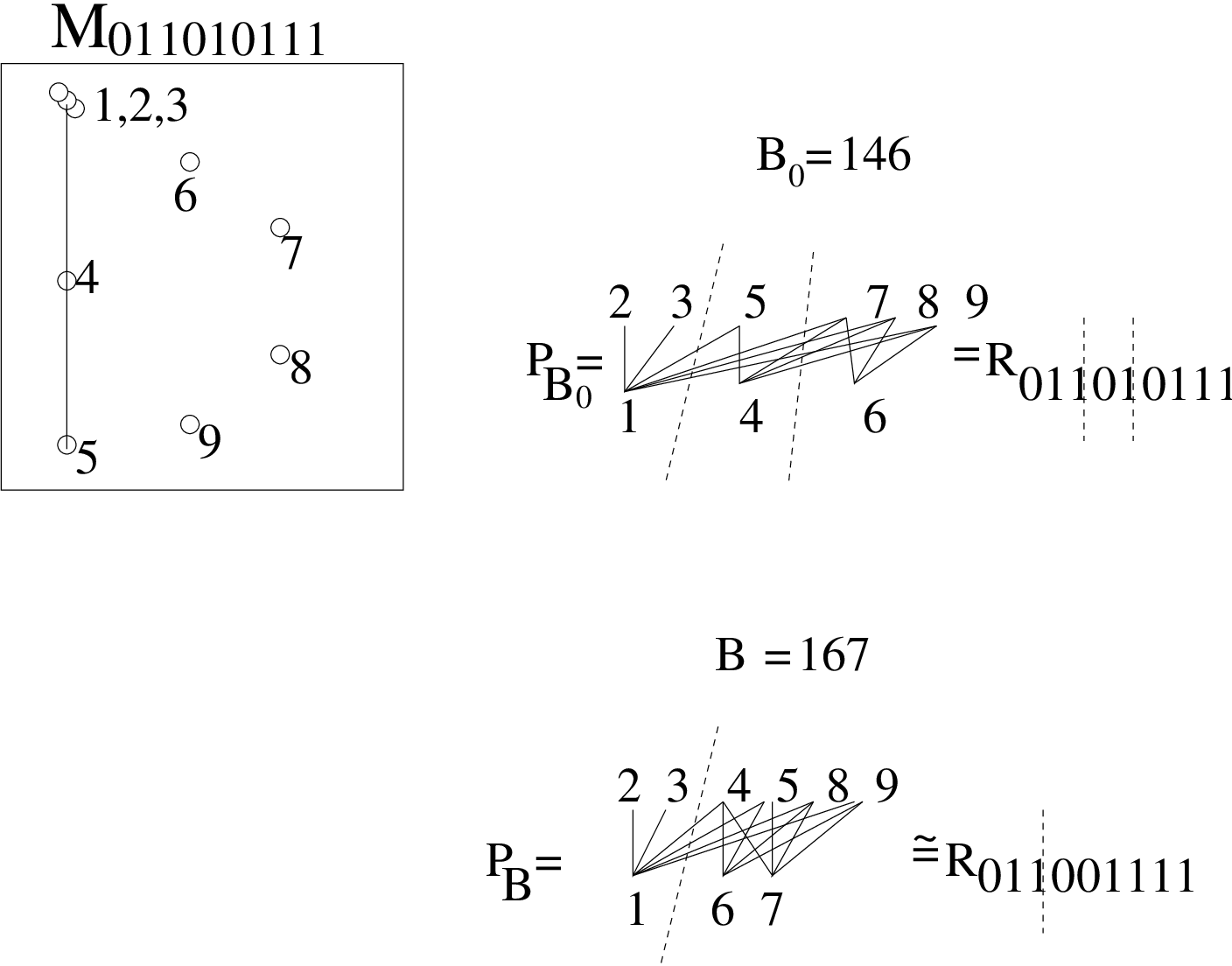}
\end{center}
 \caption{
\label{diagonal-bases-figure}
An example of the proof of Proposition~\ref{PI-triangularity}.
}
\end{figure}

\begin{proposition}
\label{PI-triangularity}
The lexicographic order $<_{lex}$ on $0\{0,1\}^{n-1}$ induced by $0 < 1$
makes the expansion of $\{ F(M_\sigma,\xx) : \sigma \in 0\{0,1\}^{n-1} \}$
in terms of the strict 
$P$-partition enumerators $\{ F(R_\sigma,\gamma_\sigma,\xx) : \sigma \in 0\{0,1\}^{n-1} \}$
triangular of the following form:
\begin{equation}
\label{Step1-expansion}
F(M_\sigma,\xx) = 
\sum_{ \tau \leq_{lex} \sigma } c_{\sigma,\tau} F(R_\tau,\gamma_\tau,\xx).
\end{equation}
where $c_{\sigma,\tau} \in \ZZ$, and $\gamma_\tau$ is any strict labelling of the poset $R_\tau$.
Furthermore, the diagonal coefficient $c_{\sigma,\sigma}$ can be expressed in terms
of the blocks $(A_1,\ldots,A_t)$ and vector $(z_1,\ldots,z_t)$ associated to $\sigma$
as follows:
$$
c_{\sigma,\sigma} = \prod_{i=1}^t \binom{|A_i|}{z_i}.
$$
\end{proposition}
\begin{proof}
We use Theorem~\ref{poset-expansion} and expand 
$$
F(M_\sigma,\xx) = \sum_{B \in \bases(M_\sigma)} F(P_B,\gamma_B,\xx).
$$
The bases $B$ of $M_\sigma$ are easily analyzed in terms of
the blocks $(A_1,\ldots,A_t)$ and vector $(z_1,\ldots,z_t)$ associated to $\sigma$
(cf. \cite[Proposition 5.1]{CrapoSchmitt1}).
Note that $M_\sigma$ will have rank $r:=z_1+\cdots+z_t$, and
it has a distinguished chain of flats
$$ 
\varnothing \subset F_1 \subset \cdots \subset F_t = [n]
$$
in which $F_i:=A_1 \sqcup A_2 \sqcup \cdots \sqcup A_i$.
Bases $B$ of $M_\sigma$ are then simply the $r$-subsets $B$ of $[n]$ that
contain for each $i=1,\ldots,t$ at most $z_1 + z_2 + \cdots + z_i$
elements from the flat $F_i$.

Given any base $B$ of $M_\sigma$, we claim that the poset $P_B$ is
isomorphic to some $R_\tau$.  To see this, we use Proposition~\ref{shifted-characterization}.
We know that $P_B$ has height at most $1$.  Relabel its minimal (resp. maximal) elements,
that is, those in $B$ (resp. $B^\ast\ast$) by
an initial (resp. final) segment of $[n]$, with those lying in block $A_i$ coming earlier than those
in block $A_j$ whenever $i < j$.  It is then easy to check that any minimal (resp. maximal)
element of $P_B$ will have its upward (resp. downward) neighbors in $P_B$ forming a
final (resp. initial) segment of $[n]$.

The diagonal terms on the right side of \eqref{Step1-expansion} come from
bases $B$ of $M_\sigma$ containing {\it exactly} $z_i$ elements of $F_i\backslash F_{i-1}$ for each $i$;
let us call these the {\it diagonal} bases of $M_\sigma$.
For example, the lexicographically earliest
base $B_0$ for $M_\sigma$ is a diagonal base, 
and it is not hard to see that $P_{B_0} = R_\sigma$ on the nose; see Figure~\ref{diagonal-bases-figure} for an example.
There are a total of $\prod_{i=1}^t \binom{|A_i|}{z_i}$ diagonal bases $B$ for $M_\sigma$,
and each has $P_B \cong P_{B_0} = R_\sigma$.

For any non-diagonal base $B$, there is some smallest index $i$ such that $B$ contains
less than $z_i$ elements of $F_i\backslash F_{i-1}$. It is not hard to see that such a $B$ will
have $P_B \cong R_\tau$ for some $\tau$ that agrees with $\sigma$ in the first
$|F_{i-1}|$ positions, that is, in the positions indexed by their first $i-1$ blocks of $\sigma$ 
(or $\tau$).
But then the $i^{th}$ block $A_i$ for $\tau$ indexes a $\{0,1\}$-substring of $\tau$
of the form $00 \cdots 011 \cdots 1$ starting with more zeroes than does 
the corresponding $i^{th}$ block $A_i$ for $\sigma$, so that $\tau <_{lex} \sigma$.
\end{proof}

\begin{example} \rm \ 
Figure~\ref{diagonal-bases-figure} illustrates the previous proof.
Here $\sigma=01101011$.  The matroid $M_\sigma$ is drawn as an affine point configuration.
Its associated chain of flats is $$F_1=123 \subset F_2=12345 \subset F_3=123456789.$$
The lexicographically first base $B_0=146$ of $M_\sigma$ is a diagonal base, having
poset $P_{B_0}$ which coincides with $R_\sigma$.
An example of a non-diagonal base $B=167$ is shown, with poset $P_B$ isomorphic to $R_\tau$ where 
where $\tau = 01101011$.  Here the smallest index $i$ for which $B$ does not contain
$z_i$ elements of $F_i - F_{i-1}$ is $i=2$, and hence $\sigma, \tau$ agree in
their first $|F_1|=3$ positions.  However the second block $A_2=456789$ in $\tau$ indexes a substring
$00111$ starting with two zeroes, while the second block $A_2=45$ in $\sigma$ indexes a substring
$01$ starting with only one zero. Hence $\tau <_{lex} \sigma$.
\end{example}

This completes Step 1 of our program:  the formula for $c_{\sigma,\sigma}$ in the previous
result only contains factors of the form $\binom{|A_i|}{z_i}$ in which $|A_i| \leq n$, so
that each of these factors divides $n!$.

\subsection{The first new basis for $\Qsym$}
\label{Q-basis-section}

In this subsection, we complete Step 2 of the proof of Theorem~\ref{refined-surjectivity-thm}
by exhibiting a new $\ZZ$-basis for $\Qsym$ that may be of independent interest.
This basis turns out to have a nice expansion property 
(Lemma~\ref{positive-triangular-expansion-property}) when one multiplies
one of its elements by $L_1=x_1+x_2+\cdots$.

This new basis comes from a family of (non-naturally, non-strictly)
labelled posets $Q_\sigma$ on $[n]$, indexed
by $\sigma$ in $0\{0,1\}^{n-1}$, which are defined recursively.  Before defining
them, we recall some standard labelled poset terminology.

\medskip

Let $P_1, P_2$ be labelled posets on label sets $A_1, A_2$ that disjointly decompose $[n]$,
that is, $[n]=A_1 \sqcup A_2$.  Their {\it disjoint sum} $P_1 + P_2$ is  
the labelled poset on label set $[n]$ keeping all order relations that were present in
$P_1$ or in $P_2$, with no new order relations between $P_1$ and $P_2$.  
Their {\it ordinal sum} $P_1 \oplus P_2$ is
obtained by from the disjoint sum by imposing further new order relations:
$p_1 < p_2$ for all $p_1 \in P_1, p_2 \in P_2$.

Now one can define the labelled posets
$Q_\sigma$ for $\sigma$ in $0\{0,1\}^{n-1}$ recursively by:
\begin{enumerate}
\item[$\bullet$] 
$Q_{\underbrace{00\cdots0}_{n\text{ zeroes}}}$ is the labelled poset on $[n]$
which is an antichain.
\item[$\bullet$] If $\sigma$ ends with a $1$, say $\sigma=\hat\sigma 1$, then
$Q_{\sigma} = Q_{\hat\sigma} \oplus (n+1)$
where $(n+1)$ is a labelled poset with one element labelled $n+1$.
\item[$\bullet$] If $\sigma$ ends with a $0$ (but is not {\it all} zeroes), say $\sigma=\hat\sigma 0$,
then $Q_\sigma$ is obtained from $Q_{\hat\sigma}$ by adding in a new element labelled $n+1$, with
only one new order relation $n < n+1$ (plus all others generated by transitivity), and then
swapping the labels of $n, n+1$.
\end{enumerate}

\begin{figure}
\begin{center}
\includegraphics[scale=0.30]{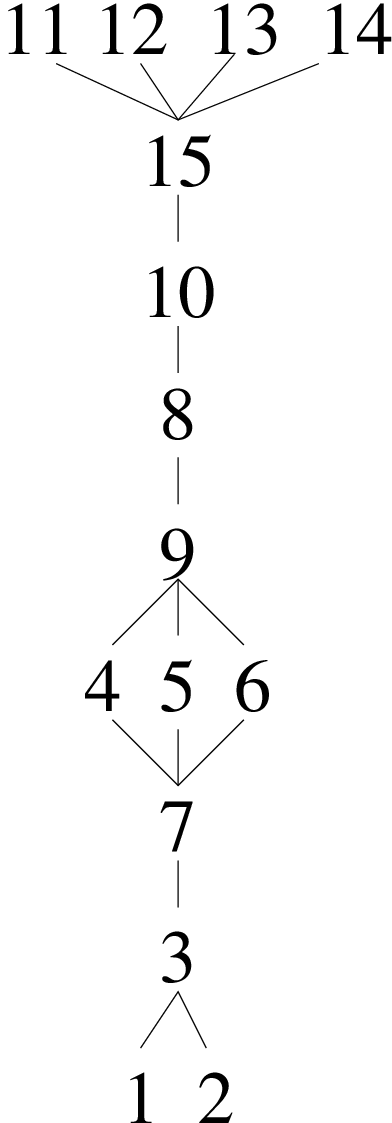}
\end{center}
 \caption{
\label{Q-poset-figure}
The poset $Q_{001100010110000}$.
}
\end{figure}

\begin{example} \rm \  
\label{Q-poset-example}
The string $\sigma$ in $0\{0,1\}^{14}$ given by
$$
\begin{matrix}
\sigma=&0&0&1&1&0&0&0&1&0&1 &1 &0 &0 &0 &0 \\
       &1&2&3&4&5&6&7&8&9&10&11&12&13&14&15
\end{matrix}
$$
has $Q_\sigma$ given by these order relations:
$$
1,2 < 3 < 7 < 4,5,6 < 9 < 8 < 10 < 15 < 11,12,13,14
$$
as depicted in Figure~\ref{Q-poset-figure}.
\end{example}

It is not hard to see that $Q_\sigma$ is always isomorphic to an 
iterated ordinal sum of a sequence of antichains.  For example, in
the poset $Q_\sigma$ of Example~\ref{Q-poset-example}, these
antichains are the induced subposets on these sets:
$$
\{1,2\}, \{3\}, \{7\}, \{4,5,6\},\{9\},\{8\},\{10\}, \{15\}, \{11, 12, 13, 14\}.
$$

\begin{remark} \rm \  
The recursive definition of $Q_\sigma$ can be rephrased, after introducing a certain
simple operation on labelled posets, which will be useful later.

For each positive integer $m$, define an operation $\psi_m$ that takes labelled posets 
on $[n]$ to labelled posets on $[n+m]$ as follows.  Given a labelled poset $P$ on
$n$, then $\psi_m(P):= P \oplus (n+m) \oplus A$ where $(n+m)$ is a labelled poset
with one element labelled $n+m$, and $A$ is an $(m-1)$-element antichain with elements labelled
$n+1, n+2,...,n+m-1$. 

To describe $Q_\sigma$ in terms of these operations, 
uniquely decompose $\sigma$ into an initial sequence of $n_0$ zeroes, and then
sequences of length $n_1,n_2,\ldots,n_p \geq 1$ of the form ${100\cdots0}$
Then 
$$
Q_\sigma:=\psi_{n_p} \cdots \psi_{n_2} \psi_{n_1} (Q_{\underbrace{00\cdots0}_{n_0\text{ zeroes}}}).
$$
\end{remark}

One then has the following proposition.

\begin{proposition}
\label{Q-basis-prop}
The $P$-partition enumerators $\{ F(Q_\sigma, \xx) : \sigma \in 0\{0,1\}^{n-1} \}$
form a $\ZZ$-basis for $\Qsym_n$.
\end{proposition}
\begin{proof}
Given $\sigma \in 0\{0,1\}^{n-1}$, let $w_\sigma$ be the linear extension of
the labelled poset $Q_\sigma$ obtained by reading each of the antichains
discussed above in the {\it reverse} of their usual numerical order.
E.g. one has  
$$
w_\sigma = 2\cdot1\quad3\quad7\cdot6\cdot5\cdot4\quad9\cdot8\quad10\quad15\cdot14\cdot13\cdot12\cdot11
$$
in the previous example, where we have indicated the positions of
descents in $w_\sigma$ by dots.

It is easily seen that
\begin{enumerate}
\item[$\bullet$] the descent set of $w_\sigma$ can be read from $\sigma$ as follows:
$$
\Des(w_\sigma)= \{i \in [n-1]:\sigma_{i+1} = 0\},
$$
and 
\item[$\bullet$] every other linear extension $w$ in $\LL(Q_\sigma)$ has
$$
\Des(w) \subsetneq \Des(w_\sigma)
$$
because at least one of the antichains discussed above must not appear in 
reverse order in $w$.  
\end{enumerate}

Hence the expansion
$$
F(Q_\sigma) = L_{\alpha(w_\sigma)} + \sum_{w \in \LL(Q_\sigma)-\{w_\sigma\}}  L_{\alpha(w)}
$$
is unitriangular with respect to the lexicographic orders on the set $0\{0,1\}^{n-1}$ and the
set of compositions $\alpha$ of $n$.
\end{proof}

\subsection{An expansion property} \ \\
It turns out that the $F(Q_\sigma, \xx)$ basis for $\Qsym$ has an interesting expansion property when one
multiplies by $L_1:=x_1 + x_2 + \cdots$.   The expansion is both nonnegative, and triangular in
a certain sense;  see Lemma~\ref{positive-triangular-expansion-property} below.

Before diving into its statement and proof, 
we introduce some notation, and observe a few simple facts about
$P$-partition enumerators.

\begin{defn} \rm \ \\
Let $P$ be a labelled poset on $n$ integers $\omega_1 <_\ZZ \ldots <_\ZZ \omega_n$.
Then the {\it standardization} $\std(P)$ of $P$ is the labelled poset on 
$[n]$ obtained by replacing the label $\omega_i$ in $P$ with
the integer $i$ for $i=1,\ldots,n$.

Given binary strings $\sigma$ and $\tau$, denote their concatenation
by $\sigma \tau$;  the most frequently used case for us will be
where $\tau=100\cdots0$ so that $\sigma\tau = \sigma100\cdots0$.
\end{defn}

The next two propositions should then be clear from Proposition~\ref{Stanley's-P-partition-result},
and will be used repeatedly without reference.

\begin{proposition}
Let $P$ be a labelled poset on $[n]$ which is an ordinal sum 
$$
P=P_1 \oplus (n) \oplus P_2
$$
in which $(n)$ is the labelled poset
with one element labelled $n$, and $P_1$ have $n_1, n_2$ elements respectively
(so that $n_1+n_2+1=n$).

Let $P'$ be the following labelled poset on $[n]$.  First form the
labelled poset $P_2'$ on $[n_2]$ obtained from $\std(P_2)$ by adding $n_1$ to all
of its labels.  Define 
$$
P' := \std(P_1) \oplus (n) \oplus P_2'.
$$
Then 
$$
F(P,\xx) = F( P',\xx) \qed.
$$

\end{proposition}

\begin{proposition}
\label{psi-map-prop}
The $\ZZ$-linear map $\psi_m: \Qsym_n \longrightarrow \Qsym_{n+m}$
defined by sending
$$
F(w,\xx)   \longmapsto F( \psi_m(w), \xx)
$$
for any permutation $w$ will also send
$$
F(Q_\sigma, \xx) \longmapsto F(Q_{\sigma 100\cdots 0},\xx)
$$
and more generally, for any labelled poset $P$ on $[n]$, sends
$$
F(P, \xx) \longmapsto F(\psi_m(P),\xx).
\qed
$$
\end{proposition}

\noindent
Note that we are slightly abusing terminology here, in using the same name $\psi_m$ for
a $\ZZ$-linear map and also for an operation on posets.

We now come to the crucial expansion property of the $F(Q_\sigma,\xx)$ basis.

\begin{lemma}
\label{positive-triangular-expansion-property}
For any $\sigma$ in $0\{0,1\}^{n-1}$,
$$
F(Q_\sigma, \xx) \cdot L_1 
= F(Q_{\sigma 0}, \xx) + \sum_{\tau <_{lex} \sigma 0} c_\tau F(Q_\tau, \xx)
\text{ with }c_\tau \in \NN.
$$
\end{lemma}

\begin{proof}
Induct on $n$. One has
\begin{equation}
\label{common-poset-expansion}
\begin{aligned}
F(Q_\sigma, \xx) \cdot L_1 &= F(Q_{\sigma} + (n+1), \xx) \\
                           &= \sum_{w \in \LL(Q_{\sigma} + (n+1))} F(w,\xx)
\end{aligned}
\end{equation}
\noindent
We analyze the set of linear extensions $\LL(Q_{\sigma} + (n+1))$.  The analysis
breaks up into two cases.

\vskip.1in
\noindent
{\bf Case 1.} $\sigma$ ends with a $1$, say $\sigma = \hat{\sigma} 1$.

In this case, $n$ is a top element of $Q_\sigma$ by construction,
and we decompose the linear extensions $w$ in $\LL(Q_{\sigma} + (n+1))$ into three sets, 
based on the location of $n+1$ relative to $n$:
\begin{enumerate}
\item[$\LL_1:$]
Those $w$ with $n+1$ occurring second-to-last, just before $n$.
\item[$\LL_2:$]
Those $w$ with $n+1$ occurring last, just after $n$.
\item[$\LL_3:$]
Those $w$ remaining, in which $n+1$ occurs at least two positions before $n$.
\end{enumerate}

\noindent
It is easy to see that $\LL_1 = \LL(Q_{\sigma 0})$.

Letting $t=(n,n+1)$ denote the transposition that swaps the labels $n, n+1$ in a labelled
poset on $[n+1]$, a little thought shows 
$$
\LL_2 \sqcup t \LL_3 = \LL( (Q_{\hat\sigma} + (n)) \oplus (n+1) ).
$$
Also, if one applies $t$ to a linear extension in which $n, n+1$ are not adjacent,
there is no effect on the descent set.  Since this is true for every linear extension
in $\LL_3$, one knows that $\LL_2 \sqcup t \LL_3$ has the same distribution of descent sets 
as $\LL_2 \sqcup \LL_3$.

Therefore, \eqref{common-poset-expansion} implies
\begin{equation}
\label{Case1-eqn-string}
\begin{aligned}
F(Q_\sigma, \xx) \cdot L_1 
                 &=F(Q_{\sigma 0}, \xx) + F( (Q_{\hat\sigma} + (n)) \oplus (n+1), \xx) \\
                 &=F(Q_{\sigma 0}, \xx) + \psi_1 ( F( Q_{\hat\sigma} , \xx) \cdot L_1 ) \\
                 &=F(Q_{\sigma 0}, \xx) + \psi_1 ( F(Q_{\hat\sigma 0},\xx) 
                    + \sum_{\hat{\tau} <_{lex} \hat\sigma 0} c_{\hat \tau} F(Q_{\hat \tau}) )\\
                 &=F(Q_{\sigma 0}, \xx) + F(Q_{\hat\sigma 01},\xx) 
                    + \sum_{\hat{\tau} <_{lex} \hat\sigma 0} c_{\hat \tau} F(Q_{\hat \tau 1}) \\
\end{aligned}
\end{equation}
where the third equality uses the inductive hypothesis.  
Since 
$$
\hat \tau 1 <_{lex} \hat \sigma 0 1 <_{lex} \hat \sigma 1 0 = \sigma 0,
$$
the last equation in \eqref{Case1-eqn-string} gives the desired conclusion.

\vskip.1in
\noindent
{\bf Case 2.} $\sigma$ ends with a $0$, say 
$\sigma = \hat{\sigma} \underbrace{100\cdots 0}_{m\text{ letters}}$.

This time we decompose the linear extensions $w$ in $\LL(Q_{\sigma} + (n+1)$ into four sets, 
again based on the location of $n+1$ relative to $n$:
\begin{enumerate}
\item[$\LL_1:$]
Those $w$ with $n+1$ at least two positions after $n$.
\item[$\LL_2:$]
Those $w$ with $n+1$ immediately after $n$.
\item[$\LL_3:$]
Those $w$ with $n+1$ immediately preceding $n$.
\item[$\LL_4:$]
Those $w$ with $n+1$ at least two positions before $n$.
\end{enumerate}

\noindent 
Note that the sets $\LL_1, \LL_4$ will have their descent set distributions unchanged
when one applies the transposition $t=(n, n+1)$ to their labels.  A little thought then shows that 
$$
\LL_3 \sqcup t\LL_1 = 
\LL(Q_{\hat{\sigma} \underbrace{100\cdots 0}_{m+1\text{ letters}}}) = \LL(Q_{\sigma 0})
$$
and
$$
\LL_2 \sqcup t \LL_4 = \LL( \psi_m( Q_{\hat\sigma} + (n)) ).
$$

Consequently \eqref{common-poset-expansion} implies
\begin{equation}
\label{Case2-eqn-string}
\begin{aligned}
F(Q_\sigma, \xx) \cdot L_1 
                 &=F(Q_{\sigma 0}, \xx) + F( \psi_m (Q_{\hat\sigma} + (n)) , \xx) \\
                 &=F(Q_{\sigma 0}, \xx) + \psi_m ( F( Q_{\hat\sigma} , \xx) \cdot L_1 ) \\
                 &=F(Q_{\sigma 0}, \xx) + \psi_m ( F(Q_{\hat\sigma 0},\xx) 
                    + \sum_{\hat{\tau} <_{lex} \hat\sigma 0} c_{\hat \tau} F(Q_{\hat \tau}) )\\
                 &=F(Q_{\sigma 0}, \xx) + F(Q_{\hat\sigma 0100\cdots 0},\xx) 
                    + \sum_{\hat{\tau} <_{lex} \hat\sigma 0} c_{\hat \tau} F(Q_{\hat \tau 100\cdots 0}) \\
\end{aligned}
\end{equation}
where the third equality uses the inductive hypothesis.  
Since
$$
\hat \tau 1 0 0 \cdots 0 <_{lex} \hat \sigma 0 1 0 0 \cdots 0 <_{lex} \hat \sigma 1 0 0 \cdots 0  = \sigma 0,
$$
the last equation in \eqref{Case2-eqn-string} gives the desired conclusion.

\end{proof}

This completes Step 2 in the proof of Theorem~\ref{refined-surjectivity-thm}.

\subsection{The second new basis for $\Qsym$}
\label{R-Q-expansion-section}

The goal of this subsection is to prove the following positive, unitriangular
expansion of the $F(R_\sigma,\xx)$ in terms of the $F(Q_\sigma,\xx)$.

\begin{theorem}
\label{R-Q-expansion}
For $\sigma$ in $0\{0,1\}^{n-1}$,
$$
F(R_\sigma,\xx) = F(Q_\sigma,\xx) + 
\sum_{ \tau <_{lex} \sigma } c_{\tau} F(Q_\tau,\xx).
$$
for some $c_{\tau}$ in $\NN$.
\end{theorem}

Note that this implies the $F(R_\sigma,\xx)$ form a $\ZZ$-basis for $\Qsym$,
which would complete Step 3 of the proof of Theorem~\ref{refined-surjectivity-thm}.

Theorem~\ref{R-Q-expansion} is simply the conjunction of assertions (i) and (ii) in the following lemma.

\begin{lemma}
For $\sigma$ in $0\{0,1\}^{n-2}$,
\begin{enumerate}
\item[(i)] 
$$
F(R_{\sigma 1},\xx) = F(Q_{\sigma 1},\xx) + 
\sum_{ \tau <_{lex} \sigma 1 } c_{\tau} F(Q_{\tau},\xx) \text{ with }c_\tau \in \NN.
$$
\item[(ii)] 
$$
F(R_{\sigma 0},\xx) = F(Q_{\sigma 0},\xx) + 
\sum_{ \tau <_{lex} \sigma 0 } c_{\tau} F(Q_{\tau},\xx) \text{ with }c_\tau \in \NN.
$$
\item[(iii)] For $\sigma$ in $0\{0,1\}^{n-m}$ and $m \geq 1$,
$$
\psi_m F(R_\sigma,\xx) = F(Q_{\sigma 100\cdots 0},\xx) + 
\sum_{ \tau <_{lex} \sigma 100\cdots 0} c_{\tau} F(Q_\tau,\xx) \text{ with }c_\tau \in \NN .
$$
\end{enumerate}
\end{lemma}
\begin{proof}
We prove all three assertions (i),(ii),(iii) by a simultaneous induction on $n$.

\vskip.1in
\noindent
{\bf Proof of (ii).}

Given $\sigma$ in $0\{0,1\}^{n-2}$, one has
$$
\begin{aligned}
F(R_{\sigma 0},\xx) 
  &= F(R_\sigma + (n), \xx) \\
  &= F(R_\sigma,\xx) \cdot L_1 \\
  &= ( F(Q_\sigma,\xx) + \sum_{\tau <_{lex} \sigma} c'_\tau F(Q_\tau,\xx) ) \cdot L_1 \\
  &= F(Q_{\sigma 0},\xx) + \sum_{\rho <_{lex} \sigma 0} c''_\rho F(Q_\rho,\xx)
      + \sum_{\tau <_{lex} \sigma} \quad 
        \sum_{\nu \leq_{lex} \tau 0} c'_\tau c''_\nu F(Q_\nu,\xx)
\end{aligned}
$$
where the third equality uses induction, and the
last equality uses Lemma~\ref{positive-triangular-expansion-property}.  
Note that the last equality implies assertion (ii).

\vskip.1in
\noindent
{\bf Proof of (iii).}

Given $\sigma$ in $0\{0,1\}^{n-m}$ and $m \geq 1$, one has
$$
F(R_\sigma,\xx) = F(Q_\sigma,\xx) + 
\sum_{ \tau <_{lex} \sigma } c_{\tau} F(Q_\tau,\xx) \text{ with }c_\tau \in \NN,
$$
by induction using assertions (i),(ii) (that is, Theorem~\ref{R-Q-expansion}).
Applying $\psi_m$ to this equality gives
$$
\begin{aligned}
\psi_m F(R_\sigma,\xx) 
&= \psi_m F(Q_\sigma,\xx) + \sum_{ \tau <_{lex} \sigma } c_{\tau} \psi_m F(Q_\tau,\xx) \\
&= F(Q_{\sigma 100 \cdots 0},\xx) + \sum_{ \tau <_{lex} \sigma } c_{\tau} F(Q_{\tau100\cdots 0},\xx)
\end{aligned}
$$
where the last equality uses Proposition~\ref{psi-map-prop}.  This gives assertion (iii).

\vskip.1in
\noindent
{\bf Proof of (i).}
Given $\sigma$ in $0\{0,1\}^{n-2}$, let
$$
J:=\{j \in [n-1]: \sigma_j=1\},
$$
so that the labelled poset $R_{\sigma 1}$ has the element labelled $n$ above
{\it all} of the elements in $[n]-J$, and above {\it none} of the element in $J$.
This means that for every linear extension $w$ in $\LL(R_{\sigma 1})$, there is 
a unique subset $I \subseteq J$ consisting of those elements appearing later
({\it i.e.} higher) in $w$ than $n$.  A little thought shows that this gives a decomposition
$$
\LL(R_{\sigma 1}) = \bigsqcup_{I \subseteq J} \LL(P^I)
$$
where $P^I:=P^I_1 \oplus (n) \oplus P^I_2$ is a labelled poset on $[n]$ 
having $P^I_1$ the restriction of $R_\sigma$ to its elements labelled by $[n-1]-I$,
and $P^I_2$ an antichain labelled by the elements of $I$.
Consequently,
$$
\begin{aligned}
F(R_{\sigma 1},\xx) 
 &= \sum_{I \subseteq J} F(P^I,\xx) \\
 &= \sum_{I \subseteq J} F(P^I_1 \oplus (n) \oplus P^I_2 ,\xx) \\
 &= \sum_{I \subseteq J} \psi_{|I|+1} F(\std(P^I_1),\xx)\\
 &= \sum_{I \subseteq J} \psi_{|I|+1} F(R_{\sigma\backslash I},\xx)
\end{aligned}
$$
where for each $I \subseteq J$, the string $\sigma \backslash I$ is obtained from the string $\sigma$
by removing all the ones that were in the positions indexed by $I$.
Hence by induction using assertion (iii) one obtains
\begin{equation}
\label{assertion-(iii)-usage}
F(R_{\sigma 1},\xx) 
 = \sum_{I \subseteq J} \left( F(Q_{(\sigma\backslash I)100\cdots 0},\xx) 
         + \sum_{\tau <_{lex} (\sigma \backslash I)100\cdots 0} c_\tau F(Q_\tau,\xx) \right)
\end{equation}
with $c_\tau$ in $\NN$.  Note that for any $I \subseteq J$ one has
$$
(\sigma \backslash I)\underbrace{100\cdots 0}_{|I|+1\text{ letters}} \leq_{lex} \sigma 1
$$
and equality occurs if and only if $I=\varnothing$.  Hence assertion (i) follows
from \eqref{assertion-(iii)-usage}.
\end{proof}

\subsection{Remarks on the bases for $\Qsym$}
We close with a few remarks on these new bases for $\Qsym$.

\begin{remark} \rm \  
Note that the $F(R_\sigma,\xx)$ basis for $\Qsym$ consists entirely of {\it naturally labelled}
$P$-partition enumerators.  This answers affirmatively the question of whether $\Qsym$ is
$\ZZ$-linearly spanned by naturally labelled $P$-partition enumerators;  note that neither
of the usual $\ZZ$-bases for $\Qsym$ (the $M_\alpha$ or $L_\alpha$) have this form.

The same question was also answered (affirmatively) in recent work of Stanley \cite{Stanley3}
who, after being queried by the authors of the current paper, produced yet another $\ZZ$-basis
for $\Qsym$ consisting of naturally labelled $P$-partition enumerators.  Given a composition
$\alpha=(\alpha_1,\ldots,\alpha_k)$ of $n$, he defined $P_\alpha$ to be the naturally labelled poset
which is the ordinal sum $A_1 \oplus \cdots \oplus A_k$, 
in which $A_i$ is an antichain on $\alpha_i$ elements
for each $i=1,2,\ldots,k$.  These posets $P_\alpha$ bear a close resemblance
to the (non-naturally) labelled posets $Q_\sigma$ defined above, in that both have simple, 
unitriangular expansions of their $P$-partition enumerators 
in terms of the $L_\alpha$-basis.  In \cite{Stanley3}, Stanley combinatorially interprets 
this upper unitriangular change-of-basis matrix between his basis and the 
$L_\alpha$-basis, as well as providing a nice (and remarkably similar) combinatorial interpretation 
for the inverse change-of-basis matrix.
\end{remark}

\begin{remark} \rm \  
The matrix $A_n$ giving the expansion of $F(R_\sigma,\xx)$ into $L_\alpha$ within $\Qsym_n$
is unimodular, and it turns out that our previous results imply a nice
$LU$-decomposition for it.  

Order the strings $\sigma$ in $0\{0,1\}^{n-1}$ with lex order,
and order the compositions $\alpha$ of $n$ also in lex order.
Then the matrix $U_n$  expanding $F(R_\sigma,\xx)$ in terms of
$F(Q_\sigma,\xx)$ will be upper unitriangular (by Theorem~\ref{R-Q-expansion}),
while the matrix $L_n$ expanding $F(Q_\sigma,\xx)$ in terms of $L_\alpha$ will be
lower unitriangular (by the proof of Proposition~\ref{Q-basis-prop}).  And $A_n=L_n U_n$.

For example, when $n=3$, this looks like
$$
A_3=\left[
\begin{matrix}
        & R_{000} & R_{001} & R_{010} & R_{011}  \\
L_{111} & 1       & 0       & 0       & 0        \\
L_{12}  & 2       & 1       & 1       & 0        \\
L_{21}  & 2       & 0        & 1       & 1       \\
L_{3}   & 1       & 1       & 1       & 1       
\end{matrix}
\right]
\quad = 
$$
$$
\left[
\begin{matrix}
        & Q_{000} & Q_{001} & Q_{010} & Q_{011} \\
L_{111} & 1       &         &         &         \\
L_{12}  & 2       & 1       &         &         \\
L_{21}  & 2       & 0       & 1       &         \\
L_{3}   & 1       & 1       & 0       & 1       
\end{matrix}
\right]
\left[
\begin{matrix}
        & R_{000} & R_{001} & R_{010} & R_{011} \\
Q_{000} & 1       & 0       & 0       & 0       \\
Q_{001} &         & 1       & 1       & 0       \\
Q_{010} &         &         & 1       & 1       \\
Q_{111} &         &         &         & 1       
\end{matrix}
\right]
$$
\end{remark}

\begin{remark} \rm \  
  We have now encountered five $\ZZ$-bases for the Hopf algebra $\Qsym$ of 
quasisymmetric functions, namely 
$$
M_\alpha, L_\alpha, F(P_\alpha,\xx), F(Q_\sigma,\xx), F(R_\sigma,\xx).
$$
Given any such basis $B_\alpha$, one might ask whether the
structure constants $c^{\alpha,\beta}_{\gamma}$ from the unique expansion
$$
B_\alpha B_\beta = \sum_{\gamma} c^{\alpha,\beta}_{\gamma} B_\gamma
$$
are always nonnegative.  For the monomial basis $M_\alpha$ and the  fundamental bases $L_\alpha$,
this property is well-known to hold and is straightforward.

Unfortunately, this property fails for the remaining three bases
$P_\alpha, R_\sigma, Q_\sigma$.  They turn out to have
some negative multiplication structure constants occurring already in (relatively) low degrees:
$$
\begin{aligned}
F(P_{(1,1)},\xx) F(P_{(1)},\xx) 
  & \left ( = F(P_{(1,1)},\xx) \cdot L_1 \right ) \\
  &= F(P_{(0,0,1)},\xx) + F(P_{(0,1,0)},\xx) - F(P_{(0,1,1)},\xx) \\ 
  & \\
F(R_{01},\xx)^2 &= 2 F(R_{0101},\xx) - F(R_{0011},\xx) \\
  & \\ 
F(Q_{010},\xx)^2 &= F(Q_{001000},\xx) + 2 F(Q_{010100},\xx) + F(Q_{001100},\xx) \\
               &\qquad + 2 F(Q_{010010},\xx) - F(Q_{001001},\xx).
\end{aligned}
$$

On the other hand, the new $\ZZ$-basis for $\Qsym$ found by Luoto which was mentioned earlier {\it does} have this property; 
see \cite[\S 4.4]{Luoto}.

\end{remark}


\end{document}